\documentclass[journal]{IEEEtran}

\usepackage{graphicx}
\usepackage{dblfloatfix} 
\usepackage{ifpdf}
\usepackage[noadjust]{cite}
\usepackage{epstopdf}

\usepackage{ifpdf}
\ifpdf
\DeclareGraphicsExtensions{.eps}
\else
\DeclareGraphicsExtensions{.eps}
\fi

\newcommand\MYhyperrefoptions{bookmarks=true,bookmarksnumbered=true,
	pdfpagemode={UseOutlines},plainpages=false,pdfpagelabels=true,
	colorlinks=true,linkcolor={black},citecolor={black},urlcolor={black},
	pdftitle={stopf-tsg-final},
	pdfsubject={},
	pdfauthor={Mohammadhafez Bazrafshan and Nikolaos Gatsis},
	pdfkeywords={}}

\ifCLASSINFOpdf
\usepackage[\MYhyperrefoptions,pdftex]{hyperref}
\else
\usepackage[\MYhyperrefoptions,breaklinks=true,dvips]{hyperref}
\usepackage{breakurl}
\fi

\ifCLASSOPTIONcompsoc
\usepackage[caption=false,font=normalsize,labelfon
t=sf,textfont=sf]{subfig}
\else
\usepackage[caption=false,font=footnotesize]{subfi
	g}
\fi

\usepackage{amsmath}
\usepackage{amsfonts} 
\ifCLASSOPTIONcompsoc
  \usepackage[caption=false,font=normalsize,labelfont=sf,textfont=sf]{subfig}
\else
  \usepackage[caption=false,font=footnotesize]{subfig}
\fi

\usepackage{url}

\usepackage[inline]{enumitem}
\newcommand{\bm}[1]{\boldsymbol{#1}}
\newcommand{\mb}[1]{\mathbf{#1}}
\newcommand{\mc}[1]{\mathcal{#1}}
\newcommand{\mbb}[1]{\mathbb{#1}}
\newcommand{\mr}[1]{\mathrm{#1}}
\newcommand{\bmat}[1]{\begin{bmatrix} #1 \end{bmatrix}}
\usepackage[mathscr]{euscript}
\usepackage{amsthm}
\theoremstyle{remark}
\newtheorem{remark}{Remark}

\usepackage{tikz}
\usetikzlibrary{dsp,fit}
\usetikzlibrary{decorations.pathmorphing} 
	\usetikzlibrary{shapes.geometric}
	\usetikzlibrary{shapes.arrows}
\usetikzlibrary{dsp,fit}
\usetikzlibrary{decorations.pathmorphing} 
\tikzstyle{block} = [draw,rectangle,thick,minimum height=3em,minimum width=5em]
\tikzstyle{branch} = [circle,inner sep=1pt,minimum size=2mm,fill=white,draw=black]
\tikzstyle{branch1} = [inner sep=1pt,minimum size=0.1mm,fill=white,draw=black]
\tikzstyle{branch1} = [circle,inner sep=0pt,minimum size=1mm,fill=black,draw=black]
\usetikzlibrary{matrix} 
\usetikzlibrary{arrows} 
\usetikzlibrary{calc} 
\tikzstyle{block} = [draw,rectangle,thick,minimum height=2em,minimum width=2em]
\tikzstyle{sum} = [draw,circle,inner sep=0mm,minimum size=2mm]
\tikzstyle{connector} = [->,thick]
\tikzstyle{rightconnector}=[<-,thick]
\tikzstyle{line} = [thick]
\tikzstyle{snakeline} = [connector, decorate, decoration={pre length=0.2cm,
                         post length=0.2cm, snake, amplitude=.4mm,
                         segment length=2mm},thick, red, ->]

\usepackage{multirow}
\usepackage{algpseudocode}
\usepackage{algorithm}
\usepackage{threeparttable} 

\begin{document}

\title{Coupling Load-Following Control with OPF}

\author{Mohammadhafez~Bazrafshan,~\IEEEmembership{Student Member,~IEEE,}%

Nikolaos~Gatsis,~\IEEEmembership{Member,~IEEE,}%
        ~Ahmad~F.~Taha,~\IEEEmembership{Member,~IEEE,}%
~Joshua~A.~Taylor,~\IEEEmembership{Member,~IEEE}%
\thanks{This material is based upon work supported by the National Science Foundation Grants ECCS-1462404 and CCF-1421583.} 
}

\maketitle

\begin{abstract}
In this paper, the optimal power flow  (OPF) problem is augmented to account for the costs associated with the load-following control of a power network.  Load-following control costs are expressed through the linear quadratic regulator (LQR). The power network is described by a set of nonlinear differential algebraic equations (DAEs).  By linearizing the DAEs   around a known equilibrium, a linearized OPF that accounts for steady-state operational constraints is formulated first. This linearized OPF is then augmented by a set of linear matrix inequalities that are algebraically equivalent to the implementation of an LQR controller.  The resulting formulation, termed LQR-OPF, is a semidefinite program which furnishes optimal steady-state setpoints and an optimal feedback law to steer the system to the new steady state with minimum load-following control costs. Numerical tests demonstrate that the setpoints computed by LQR-OPF result in lower overall costs and frequency deviations compared to the setpoints of a scheme where OPF and load-following control are considered separately.
\end{abstract}

\IEEEpeerreviewmaketitle

\section{Introduction}

\IEEEPARstart{C}{apacity} expansion and generation planning,  economic dispatch, frequency regulation, and automatic generation control (AGC) are decision making problems in power networks that are solved over different time horizons. These problems range from decades in planning to several seconds in transient  control. Although decisions made in shorter time periods may negatively affect performance over longer time periods, these problems  have traditionally been treated separately. For example,  optimization is used for optimal power flow (OPF) and economic dispatch while feedback control theory is used for frequency regulation. 

This paper aims to integrate two crucial power network problems with different time-scales. The first problem is the steady-state OPF, whose decisions are updated every few minutes (e.g., 5 minutes for real-time market balancing). The second is the problem of load-following control that spans the time-scale of several seconds to one minute.  Load-following control, also known as secondary frequency control \cite{GomezExposito2009}, is responsible for maintaining system frequencies at  nominal values during normal load fluctuations \cite{Galiana2005}.

For a forecasted load level, optimal steady-state setpoints that also minimize costs of load-following control are sought.  Steady-state costs account for  generator power outputs. 
Control costs account for the action required to drive the deviation of frequency and voltage signals from their optimal OPF setpoints to zero via the linear quadratic regulator (LQR). 
 The proposed formulation also provides, as an output,  a feedback  law to guide the system dynamics to optimal OPF setpoints.

\subsection{Literature Review}
Recent research efforts, organized below in two categories,  have showcased the economic and technical merits of jointly tackling steady-state and control problems in power networks.  

The first category focuses on frequency regulation with a view towards the  economic dispatch  \cite{AndreassonDimarogonasSandbergJohansson2014, ShafieeGuerreroVasquez2014, DorflerSimpsonBullo2016,DorflerGrammatico2016,
ZhaoTopcuLiLow2014, LiZhaoChen2016, 
ZhangPapachristodoulou2015, MalladaZhaoLow2015}. These works design controllers based on feedback from frequency measurements and guarantee the stability of system dynamics---characterized by the swing equation---while ensuring that  system states  converge to steady-state values  optimal for some form of the economic dispatch problem.   

The designs of these feedback control laws may come from averaging-based controllers as in \cite{AndreassonDimarogonasSandbergJohansson2014, ShafieeGuerreroVasquez2014, DorflerSimpsonBullo2016}.  By leveraging a continuous-time version of a dual algorithm, generalized versions of the aforementioned controllers are developed in~\cite{DorflerGrammatico2016},  including the  classical AGC as a special case.  Similar  control designs come from interpretation of network dynamics as iterations of a primal-dual algorithm that solves an OPF. For example,~\cite{ZhaoTopcuLiLow2014} introduces a  primary frequency control that minimizes a reverse-engineered load disutility function, and~\cite{LiZhaoChen2016}  introduces a modified AGC that solves an economic dispatch with limited operational constraints. 

A frequency control law that further solves economic dispatch with nonlinear power flows in tree networks is devised in \cite{ZhangPapachristodoulou2015} by 
borrowing virtual dynamics from the KKT conditions.   A novel method to tackle a load-side control problem including linear equality or inequality constraints for arbitrary network topologies is presented in~\cite{MalladaZhaoLow2015}.  The chief attractive feature of \cite{AndreassonDimarogonasSandbergJohansson2014, ShafieeGuerreroVasquez2014, DorflerSimpsonBullo2016,DorflerGrammatico2016,
	ZhaoTopcuLiLow2014, LiZhaoChen2016, 
	ZhangPapachristodoulou2015, MalladaZhaoLow2015}  is that they afford decentralized or distributed implementations. Moreover, works such as \cite{DorflerGrammatico2016} and \cite{ZhangPapachristodoulou2015} develop controllers accounting for nonlinear power flow equations. 

The second category focuses on OPF variations with enhanced stability measures.  The goal is to obtain optimal steady-state setpoints less vulnerable to disturbances, rather than seeking a stabilizing control law.  This goal is achieved by  incorporating  additional constraints in OPF that account for the stability of the  steady-state optimal point.

In the context of transient stability, a classical reference  is \cite{GanThomasZimmerman2000} where  system differential equations are converted to algebraic ones and added to the OPF.  More recently,   trajectory sensitivity analysis has been extensively used to accommodate transient stability specifications in the form of linear constraints within the OPF.  Typically, these constraints are based on stability margins obtained from the (extended) equal-area criterion  \cite{TangSun2016, Xu2016,Xu2017,XinXu2017}.  In particular, \cite{TangSun2016} augments the OPF with automated computationally-efficient rotor angle constraints around a base stable trajectory. Transient-stability constrained OPF formulations that are robust to uncertainty in load dynamics and wind power generation have also been developed \cite{Xu2016,Xu2017}.  Load-shedding minimization has been pursued in \cite{XinXu2017}, where in addition to angle stability margins, trajectory sensitivities have been used to approximate constraints on voltage and frequency security margins. 

In the realm of small-signal stability,  distance from rotor instability is guaranteed in \cite{Minano2011} by providing a set of stressed load conditions as a supplement to  OPF.    The spectral abscissa, that is, the largest real part of  system state matrix eigenvalues, is upper bounded by a negative number  in \cite{LiQiWangWeiBaiQiu2016}, yielding a non-smooth OPF.  A simpler optimization problem is pursued in~\cite{MalladaTang2013}  using the pseudo-spectral abscissa as stability measure. Based on  Lyapunov's stability theorem and the system state matrix, \cite{Yang2013}  incorporates small-signal stability constraints into the OPF. Since the state matrix is a function of the steady-state variables, the overall formulation becomes a nonlinear and nonconvex semidefinite program (SDP).

\subsection{Paper contributions and organization}
The paper contributions are as follows:
\begin{itemize}[leftmargin=*]
	\item An OPF framework is proposed that solves for optimal steady-state  setpoints and also provides an optimal load-following control law to drive the system to those  setpoints.  Optimality of the control law is appraised by an integral cost on time-varying deviation of system states and controls  from their  optimal setpoints. 	An  LQR controller is then applied which minimizes this cost by providing a feedback law that is a  linear combination of system state deviations.  LQR  has previously been used in the context of megawatt-frequency control \cite{FoshaElgerd1970LQR} and control of oscillatory dynamics \cite{SinghPal2016}. The proposed framework,  in distinction,  accounts for   LQR costs from within the steady-state  time-scale. 		
	\item The setpoints computed by the proposed formulation, termed LQR-OPF, result in lower overall costs and frequency fluctuations compared to the setpoints of a scheme where OPF and load-following control are decoupled. The proposed framework also allows  time-varying control costs to  be dependent on steady-state variables.  This dependence enables subsequent regulation pricing schemes similar to \cite{TaylorNayyarCallawayPoolla2013} where for example, the cost of frequency regulation is made dependent on steady-state power generation. 
	\item In comparison to \cite{AndreassonDimarogonasSandbergJohansson2014, ShafieeGuerreroVasquez2014, DorflerSimpsonBullo2016,DorflerGrammatico2016,
		ZhaoTopcuLiLow2014, LiZhaoChen2016, 
		ZhangPapachristodoulou2015, MalladaZhaoLow2015},  the formulation in this paper includes voltage dynamics, reactive powers,  AC power flows, and a more realistic dynamical model of the synchronous generator that distinguishes between generator internal and external quantities. Costs incurred due to deviations of system states and controls from their optimal setpoints are  also accounted for.	In relation to \cite{GanThomasZimmerman2000, TangSun2016, Xu2016,Xu2017,XinXu2017, Minano2011, MalladaTang2013, LiQiWangWeiBaiQiu2016, Yang2013}, the proposed approach incorporates load-following control constraints into the OPF; but as an extra output,  the required control law to steer the system to stability is also provided.
\item The effectiveness of the derived control law is demonstrated via numerical simulations on the power system described by \textit{nonlinear differential-algebraic equations} (DAEs). We further exhibit that the steady-state setpoints provided by the proposed LQR-OPF yield smaller overall system costs even when the traditional AGC is used for load following.
\end{itemize}

In this paper, in order to  obtain a stabilizing feedback law,  nonlinear DAEs are linearized  around a known equilibrium point  with respect to the system states and algebraic variables.  It should be emphasized that this linearization is different than that of trajectory sensitivity analysis. In the latter, the system dynamics are linearized around a known trajectory with respect to system initial conditions and parameters~\cite{HiskensPai2002,TangMcCalley2013}. Formulations in~\cite{TangSun2016,Xu2016,Xu2017} that use trajectory sensitivity analysis are suitable for transient stability where it is required to limit generators' angle separation subject to contingencies and large disturbances, e.g., line trips due to three-phase faults.  On the other hand, the load-following or secondary control, which is considered in this paper, aims at maintaining frequencies at their nominal value during normal load changes.

The paper is organized as follows. The  power system model  is laid out in Section \ref{sec:network}, followed by a description of a generalized OPF.  System linearization is pursued in Section \ref{sec:systemlin}.  The proposed formulation coupling OPF and load-following control is detailed in Section \ref{sec:coupledform}.  Specific generator models,  power flow equations, and connections to the standard OPF are provided in Section \ref{sec:genmodelpfopf}.  Section \ref{sec:numtests} numerically verifies the merits of the proposed method.  Section~\ref{sec:conclusion} provides pointers for integrating more  power system applications into our proposed framework as future work.

\section{Power System Model}
\label{sec:network}
Consider a power network with $N$ buses where $\mathcal{N}  := \{1,\ldots,N\}$ is the set of nodes.  Define the partition $\mathcal{N}=\mathcal{G} \cup \mathcal{L}$ where $\mathcal{G}=\{1, \ldots, G\} $  collects  $G$  buses that contain generators (and possibly also loads) and $\mathcal{L}=\{G+1, \ldots, G+L\}$ collects the remaining $L$  load-only buses. Notice that $N= G+L$. 
For a generator $i \in \mc{G}$ with $n_s$ states and $n_c$ control inputs, denote by  $\mb{x}_i(t) \in \mbb{R}^{n_{s}}$ the time-varying vector of state variables, and denote by $\mb{u}_i(t) \in \mbb{R}^{n_{c}}$ the   time-varying control inputs. For example, adopting a fourth-order model yields $n_{s}=4$ and $n_{c}=2$ with $\mb{x}_i(t)=\{\delta_i(t),\omega_i(t), e_i(t), m_i(t)\}$ and $\mb{u}_i(t)=\{r_i(t), f_i(t)\}$ where $\delta_i(t)$, $\omega_i(t)$, $e_i(t)$, $m_i(t)$, $r_i(t)$, and $f_i(t)$ respectively denote the generator internal phase angle, rotor electrical velocity,  internal EMF,  mechanical power input, reference power setting,  and internal field voltage. Further details are given in Section~\ref{subsec:genmodel}.

Denote by $\mb{a}_i(t)$ the vector of algebraic variables. For load nodes $i \in \mc{L}$, $\mb{a}_i(t)=\{v_i(t), \theta_i(t)\}$, where $v_i(t)$ and $\theta_i(t)$ denote the terminal load voltage and phase angle. 
For generator nodes  $i \in \mc{G}$, $\mb{a}_i(t)=\{p_{g_i}(t), q_{g_i}(t), v_i(t), \theta_i(t)\}$, where $p_{g_i}(t)$, $q_{g_i}(t)$,  $v_i(t)$, and $\theta_i(t)$ respectively denote generator real and reactive power, terminal voltage and phase angle. For brevity, the dependency of variables $\mb{x}_i$, $\mb{a}_i$, and $\mb{u}_i$ on $t$ is dropped, and the notations $\mb{x}:=\{\mb{x}_i\}_{i \in \mc{G}} \in \mbb{R}^{n_sG}$, $\mb{u}:=\{\mb{u}_i\}_{i \in \mc{G}} \in \mbb{R}^{n_cG} $, and $\mb{a}:=\{\mb{a}_i\}_{i\in \mc{N}} \in \mbb{R}^{2N+2G}$ are introduced. Finally, let $\mb{z}=\{\mb{x}, \mb{a}, \mb{u}\} \in \mbb{R}^{(n_s+n_c+2)G+2N}$.  
The dynamics  of a power system can be captured by a set of nonlinear  DAEs
\begin{subequations}
\label{eqngroup:DAEs}
\begin{IEEEeqnarray}{rCl}
\dot{\mb {x}} &=& \mb {g}(\mb {x}, \mb {a}, \mb {u}), \label{eqn:genericdifferential}\\
\mb{d}& =& \mb {h} ( \mb{x}, \mb {a}), \label{eqn:genericalgebraic}
\end{IEEEeqnarray}
\end{subequations}
where $\mb{g}:\mbb{R}^{(n_s+n_c+2)G+2N} \rightarrow \mbb{R}^{n_sG}$ is given by adopting an appropriate dynamical model of the generator,  and $\mb{h}: \mbb{R}^{(n_s+2)G+2N} \rightarrow \mbb{R}^{2G+2N}$ includes generator algebraic equations as well as the network power flow equations.  Vector $\mb{d}$ collects all the network loads as well as leading zero entries coming from two generator algebraic equations per generator. A particular example  of the mapping $\mb{g}$ is provided in Section \ref{subsec:genmodel}; for the corresponding form of the mapping $\mb{h}$ and the vector $\mb{d}$ see Section \ref{subsec:powerflowequations}.

Given steady-state load conditions $\mb {d}^{\mr{eq}}$, the system steady-state operating point  is represented by an equilibrium of the DAEs~\eqref{eqngroup:DAEs}.
By setting $\dot{\mb{x}}=\mb{0}$ and allowing $\mb{x}$, $\mb{a}$, and $\mb{u}$ to reach steady states $\mb{x}^\mr{eq}$, $\mb{a}^\mr{eq}$, and $\mb{u}^\mr{eq}$, a system of $(n_s+2)G+2N$  algebraic  equations  in $(n_s+n_c+2)G+2N$ variables is derived:\begin{subequations}
\label{eqngroup:loadflow1}
\begin{IEEEeqnarray}{rCl}
\mb {0}&=& \mb {g}( \mb {x}^\mr{eq}, \mb {a}^\mr{eq}, \mb {u}^{\mr{eq}}),\\
\mb {d}^{\mr{eq}}&=& \mb {h}( \mb {x}^{\mr{eq}}, \mb {a}^{\mr{eq}}). \IEEEeqnarraynumspace \label{eqn:loadflow1alg} 
\end{IEEEeqnarray}
\end{subequations}

Let $\mc{F}(\mb {d}^{\mr{eq}})$ denote the set of solutions to \eqref{eqngroup:loadflow1}, where the dependency on the load conditions is made explicit. 
Suppose now that  $\mb{x}^\mr{eq}, \mb{a}^\mr{eq}, \mb{u}^\mr{eq}$ are to be jointly optimized so that a certain objective function $c(\mb{x}^\mr{eq}, \mb{a}^\mr{eq}, \mb{u}^\mr{eq})$  is minimized. This leads to a generalized OPF [for clarity,  the notation $(\mb{x}^s, \mb{a}^s, \mb{u}^s)$ is used to denote optimization variables, and $(\mb{x}^\mr{eq}, \mb{a}^\mr{eq}, \mb{u}^\mr{eq})$ is used to generically denote a DAE equilibrium]. 
\begin{subequations}
\label{eqngroup:generalopf}
\begin{IEEEeqnarray}{rCl}
\IEEEeqnarraymulticol{3}{l}{\min_{\mb{x}^{s}, \mb{a}^s, \mb{u}^s}~~ c\left( \mb{x}^s, \mb{a}^s, \mb{u}^s \right)} \label{eqn:generalopfobj} \IEEEeqnarraynumspace \\
\text{subj. to} ~~~ \mb{0} &=& \mb{g}( \mb{x}^s, \mb{a}^s, \mb{u}^s),  \label{eqn:generalopfg} \IEEEeqnarraynumspace \\
 \mb{d}^s &=& \mb{h}( \mb{x}^s, \mb{a}^s),  \label{eqn:generalopfh} \IEEEeqnarraynumspace \\
 \mb{a}^s &\in& \mc{A}, \IEEEeqnarraynumspace \label{eqn:generalopfextraconstraints}
\end{IEEEeqnarray}
\end{subequations}
where \eqref{eqn:generalopfextraconstraints} are the algebraic variable constraints on voltage magnitudes, line flow limits, and line current capacities. Note that the parameter vector $\mb{d}^s$ (which includes the constant-power loads) is the input to \eqref{eqngroup:generalopf}. The term \emph{generalized} refers to the fact that \eqref{eqngroup:generalopf} considers models of generators within the OPF, see e.g., \cite{Molzahn_2017} for a recent OPF example with machine models.  The connection between \eqref{eqngroup:generalopf} and the standard OPF is explained in Section~\ref{subsec:opf}. 
 The  OPF problem \eqref{eqngroup:generalopf} guarantees optimal steady-state operating costs, but does not provide minimal control costs.  Prior to introducing a  formulation that bridges stability with OPF,  linear approximation of the system dynamics is required which is presented in the next section.

\section{Linear Approximation of System Dynamics}
\label{sec:systemlin}
To obtain an approximate dynamic,  \eqref{eqngroup:DAEs} is linearized around a known operating point $\mb{z}^0:=(\mb{x}^0, \mb{a}^0, \mb{u}^0) \in \mc{F}(\mb{d}^0)$. For example, the point $\mb{z}^0$   can be a solution of the load-flow corresponding to an operating point known to the system operator.  The motivation behind this selection is to obtain tractable constraints to augment \eqref{eqngroup:generalopf} as it allows the stability constraints to take the form of properly formulated linear matrix inequalities (LMIs).  The derived control law is eventually applied to the nonlinear DAEs \eqref{eqngroup:DAEs}, rather than the linearization derived in this section.

Consider a generic \emph{equilibrium} point $\mb{z}^0 \in \mc{F}(\mb{d}^0)$, where $\mb{d}^0$ is a known load vector.  That is, the following holds:
\begin{subequations}
	\label{eqngroup:zeroequilibrium}
 \begin{IEEEeqnarray}{rCl}
\mb{0}& =& \mb{g}(\mb{x}^0, \mb{a}^0, \mb{u}^0), \label{eqn:gzeroequilibrium} \\
\mb{d}^0 &=& \mb{h}(\mb{x}^0, \mb{a}^0). \label{eqn:hzeroequilibrium}
\end{IEEEeqnarray}
\end{subequations}

Define $\Delta \mb{d}^s:=  \mb{d}^{s}- \mb{d}^0$ as the step-change difference between $\mb{d}^{s}$, the load for which the generalized OPF \eqref{eqngroup:generalopf} is to be solved, and $\mb{d}^0$, the generic load that will be used for linearization.
Equations ~\eqref{eqn:genericdifferential} and \eqref{eqn:genericalgebraic} can be linearized  around $(\mb{x}^0, \mb{a}^0, \mb{u}^0)$ by setting $(\mb{x}, \mb{a}, \mb{u}) = (\mb {x}^0, \mb{a}^0, \mb{u}^0)+ (\Delta \mb{x}, \Delta \mb{a}, \Delta \mb{u})$:
\begin{subequations}
\label{eqngroup:genericlinear}
\begin{IEEEeqnarray}{rCl}
\Delta \dot{\mb{x}} &=& \mb{g}_{\mb{x}} (\mb{z}^0) \Delta \mb{x} + \mb{g}_{\mb{a}}(\mb {z}^0)  \Delta \mb{a} + \mb{g}_{\mb{u}}(\mb{z}^0)  \Delta \mb{u}, \label{eqn:genericlineardiff} \\
\Delta \mb{d}^s&=& \mb{h}_{\mb{x}}(\mb{x}^0, \mb{a}^0)  \Delta \mb{x} + \mb{h}_{\mb{a}}(\mb{x}^0, \mb{a}^0)  \Delta \mb{a}, \label{eqn:genericlinearalg} 
\end{IEEEeqnarray}
\end{subequations}
where the notation $\mb{g}_{\mb{x}}$ defines the Jacobian with respect to $\mb{x}$, and $\mb{g}_{\mb{a}}$, $\mb{g}_{\mb{u}}$, $\mb{h}_{\mb{x}}$,  and   $\mb{h}_{\mb{a}}$ are similarly defined.  It is also understood that $\Delta \mb{x}, \Delta \mb {a}, \Delta \mb{u}$ are functions of time. 

Equation~\eqref{eqngroup:genericlinear} represents a set of linear DAEs. Next, \eqref{eqngroup:genericlinear} is leveraged in order to  develop a proper linear dynamical system that represents the network without algebraic constraints. In particular, assuming invertibility of $\mb{h}_{\mb{a}}(\mb{x}^0, \mb{a}^0)$,   \eqref{eqn:genericlinearalg} can be solved for  the algebraic variables as
\begin{IEEEeqnarray}{rCl}
\Delta \mb {a} =  - \mb{h}_{\mb{a}}^{-1}(\mb{x}^0, \mb{a}^0) \left(-\Delta \mb{d}^s+ \mb{h}_{\mb{x}} (\mb{x}^0, \mb{a}^0)   \Delta \mb{x} \right). \IEEEeqnarraynumspace \label{eqn:eliminatealgebraiclinear}
\end{IEEEeqnarray}
The assumption on invertibility of $\mb{h}_{\mb{a}}(\mb{x}^0, \mb{a}^0)$ is very mild in the sense that it holds for practical networks and for various operating points; see also \cite{MalladaTang2013} and references therein for sufficient conditions in a similar construction. Then,  substituting \eqref{eqn:eliminatealgebraiclinear} in \eqref{eqn:genericlineardiff} yields
\begin{IEEEeqnarray}{rCl}
\Delta \dot{\mb {x}} & = & \mb{A}(\mb{z}^0) \Delta \mb{x} + \mb{B}(\mb{z}^0 )\Delta \mb {u} + \mb {g}_{\mb {a}}(\mb {z}^0) \mb {h}_{\mb {a}}^{-1} (\mb {x}^0, \mb{a}^0) \Delta \mb {d}^s, \label{eqn:newlineardynamical} \IEEEeqnarraynumspace
\end{IEEEeqnarray}
where
$\mb {A}(\mb {z}^0)=\mb {g}_{\mb {x}} (\mb {z}^0)-\mb {g}_{\mb {a}}(\mb {z}^0)\mb {h}_{\mb {a}}^{-1}(\mb {x}^0, \mb{a}^0)\mb {h}_{\mb {x}}(\mb {x}^0, \mb{a}^0),$ and $ \mb {B}(\mb {z}^0)= \mb {g}_{\mb {u}}(\mb {z}^0)$.

\section{Coupling Load-Following Control with OPF}
\label{sec:coupledform}
In this section, an optimal control problem using the LQR  is presented to stabilize the linear dynamical system \eqref{eqn:newlineardynamical} to an equilibrium point, while ensuring minimal steady-state operating costs for the equilibrium of the linearized dynamics.  The equilibrium of the linearized dynamics is denoted by $\mb{z}^s = (\mb{x}^s, \mb{a}^s, \mb{u}^s)$ to distinguish it from the equilibrium of the true nonlinear system in e.g., \eqref{eqngroup:loadflow1}.    Section \ref{sec:linear-approx} leverages  the linearized dynamical system~\eqref{eqn:newlineardynamical}  to obtain a \emph{linear approximation} of the constraints~\eqref{eqn:generalopfg} and~\eqref{eqn:generalopfh}---this yields a linearized OPF.
The objective of Section \ref{sec:lqr-opf} is to incorporate stability measures with respect to the dynamical system~\eqref{eqn:newlineardynamical} in the linearized OPF. The resulting formulation outputs optimal steady-state values $\mb{z}^{s}=(\mb{x}^s, \mb{a}^s, \mb{u}^s)$ together with the optimal control law that drives the dynamic variables to stability.

\subsection{Linear approximation of the generalized OPF}
\label{sec:linear-approx}
The steady state of~\eqref{eqngroup:genericlinear} obtained by setting  $\Delta \dot{\mb {x}}=\mb {0}$ yields the following  system:
\begin{subequations}
\label{eqngroup:linearss}
\begin{IEEEeqnarray}{rCl}
\mb {0} &=& \left[\mb {g}_{\mb {x}}(\mb {z}^0) , \mb {g}_{\mb {a}} (\mb {z}^0), \mb {g}_{\mb {u}} (\mb {z}^0) \right] (\mb {z}^s - \mb {z}^0), \label{eqn:linearssg}\\ 
 \Delta \mb {d}^s &=& \left[\mb {h}_{\mb {x}} (\mb {x}^0, \mb{a}^0), \mb {h}_{\mb {a}} (\mb {x}^0, \mb{a}^0)\right] \left( (\mb {x}^s, \mb{a}^s) - (\mb{x}^0, \mb{a}^0)\right)  \label{eqn:linearssh}  \IEEEeqnarraynumspace
\end{IEEEeqnarray}
\end{subequations}
The equations in \eqref{eqngroup:linearss} are therefore a linear approximation of constraints \eqref{eqn:generalopfg} and \eqref{eqn:generalopfh}, which leads to a linearized version of OPF, formulated as follows:
\begin{IEEEeqnarray}{lrCl}\label{eqngroup:linearizedopf}
\hspace{-1.1cm}\textbf{Linearized OPF}: \notag &&&\\ \: \min_{\mb{x}^s, \mb{a}^s, \mb{u}^s} c\left( \mb{x}^s, \mb{a}^s, \mb{u}^s \right) \: \text{subj. to}\:\, \eqref{eqngroup:linearss} \:\, \text{and} \: \, \mb{a}^s \in \mc{A}. 
\end{IEEEeqnarray}
Notice that \eqref{eqngroup:linearizedopf} has the linearized version of power flow equations as part of its constraint set.

\subsection{LQR-OPF formulation}
\label{sec:lqr-opf}  
The previous section derived~\eqref{eqngroup:linearizedopf}, which is the linear approximation of the generalized OPF \eqref{eqngroup:generalopf}. Likewise, the linear system in \eqref{eqn:newlineardynamical} is the linear approximation of the system dynamics~\eqref{eqn:genericdifferential} and~\eqref{eqn:genericalgebraic}  around $\mb {z}^0$.  This section augments the linearized OPF \eqref{eqngroup:linearizedopf} with optimal control of the dynamical system in \eqref{eqn:newlineardynamical}.  

Writing \eqref{eqn:newlineardynamical} at its equilibrium $\mb{z}^s$, that is,  setting $\Delta\dot{\mb{x}}=\mb{0}$, $\Delta \mb{x}=\Delta \mb{x}^s$, and $\Delta \mb{u}=\Delta \mb{u}^s$, yields
\begin{IEEEeqnarray}{rCl}
\mb {0} = \mb A(\mb {z}^0) \Delta \mb {x}^s + \mb B(\mb {z}^0) \Delta \mb {u}^s + \mb {g}_{\mb {a}}(\mb {z}^0) \mb {h}_{\mb {a}}(\mb {x}^0, \mb{a}^0)^{-1} \Delta \mb {d}^s \label{eqn:dynamicsinfty}. \IEEEeqnarraynumspace
\end{IEEEeqnarray}
Subtracting \eqref{eqn:dynamicsinfty} from \eqref{eqn:newlineardynamical} yields the following system with the new state variable $\Delta \mb {x}':=\Delta \mb {x} - \Delta \mb {x}^s$ and new control variable $\Delta \mb {u}'= \Delta \mb {u} - \Delta \mb {u}^s$:
\begin{IEEEeqnarray}{rCl}
\Delta \dot{\mb {x}}' = \mb A(\mb {z}^0) \Delta \mb {x}' + \mb B(\mb {z}^0) \Delta \mb {u}'.  \label{eqn:newlineardynamicalprime}
\end{IEEEeqnarray}
The initial conditions of \eqref{eqn:newlineardynamicalprime} are  $\Delta \mb{x}'(0) = \Delta \mb{x} (0) - \Delta \mb{x}^s =\mb{x} (0) - \mb{x}^s $. When $\Delta \mb{x}$ in~\eqref{eqngroup:genericlinear} or~\eqref{eqn:newlineardynamical} converges to $\Delta \mb{x}^s$, $\Delta \mb{x}'$ in~\eqref{eqn:newlineardynamicalprime} converges to $\mb{0}$. The objective is to penalize the control effort to drive $\Delta \mb{x}'$ to $\mb{0}$.

To this end, weights on the control actions and state deviations are considered. In particular, an optimal LQR controller to drive the system states $\Delta \mb {x}'$ through a state-feedback control $\Delta \mb {u}'$ to their zero steady-state values is considered. By augmenting the linearized  OPF in \eqref{eqngroup:linearizedopf}, this LQR-based OPF  is written as follows:
\begin{subequations}
\label{eqngroup:LQR1}
\begin{IEEEeqnarray}{rCl}
\min_{\substack{\mb{z}^s=\{\mb x^s, \mb a^s, \mb u^s\}\\ \Delta \mb {x}',\Delta \mb {u}'}}  \hspace{-0.2cm}c(\mb {z}^s)  && +  \frac{T_{\mr{lqr}}}{2}  \int\limits_{0}^{t_f} \Delta \mb {x'}^{\top} \mb{Q} \Delta \mb {x}' + \Delta \mb {u'}^{\top} \mb{R} \Delta \mb {u}'dt \IEEEeqnarraynumspace \\
\text{subj.~to} \quad  \: && \Delta \dot{\mb {x}}' = \mb {A}(\mb {z}^0) \Delta \mb {x}' + \mb{B}(\mb{z}^0) \Delta \mb {u}' \IEEEeqnarraynumspace \\
&& \Delta \mb {x}'(0)= \mb{x}(0) - \mb{x}^s \label{eqn:lqr1deltaxprime0}\\
 & & \eqref{eqn:linearssg},\eqref{eqn:linearssh},  \: \mb {a}^s \in \mc{A}, \IEEEeqnarraynumspace \label{eqn:lqr1extra}
\end{IEEEeqnarray}
\end{subequations}
where $\mb{Q}$ and $\mb{R}$ are positive definite matrices penalizing state and control actions deviations; $T_{\mr{lqr}}$ is a scaling factor to compensate for the time-scale of the stability control problem; $t_f$ is the optimization time-scale for the OPF problem which is typically in minutes. Since $t_f$ is in minutes, the finite horizon LQR problem can be replaced with an infinite horizon LQR formulation (i.e., $t_f=\infty$) as the solution to the Riccati equation reaches steady state~\cite{boydLMI}. In \eqref{eqn:lqr1deltaxprime0}, if  the system is  assumed to operate at $\mb{x}^0$, then $\Delta \mb{x}(0) = \mb{0}$, yielding $\Delta \mb{x}'(0) = \mb{x}^0 - \mb{x}^s$.

Regulating state and control actions can be dependent on the operating points to encourage smaller steady-state variations. For example, regulating a generator's frequency becomes more costly as the real power generation increases~\cite{TaylorNayyarCallawayPoolla2013}. This coupling is captured by considering $\mb{Q}^{-1}$ and $\mb{R}^{-1}$ to be affine functions of steady-state variable $\mb{z}^s$. In particular, assume $\mb{Q}^{-1}$ and $\mb{R}^{-1}$ to be diagonal matrices as follows
\begin{subequations}
\label{eqngroup:generalqrinv}
\begin{IEEEeqnarray}{rCl}
\mb{Q}(\mb{z}^s)^{-1}= \mb{Q}_0^{-1}+\mb{Q}_1^{-1} \mr{diag}(\mb{a}^s) , \label{eqn:generalqinv} \\
\mb{R}(\mb{z}^s)^{-1} = \mb{R}_0^{-1} + \mb{R}_1^{-1} \mr{diag}(\mb{a}^s). \label{eqn:generalrinv}
\end{IEEEeqnarray}
\end{subequations}
Matrices $\mb{Q}_0$, $\mb{Q}_1$, $\mb{R}_0$, and $\mb{R}_1$ are selected so that $\mb{Q}(\mb{z}^s)$ and $\mb{R}(\mb{z}^s)$ are positive definite over the bounded range of $\mb{a}^s$ prescribed by the set $\mc{A}$. 

The corresponding infinite horizon LQR augmenting the linearized OPF \eqref{eqngroup:LQR1} can be written as an SDP, as follows:\begin{subequations}
\label{eqngroup:lqropf}
\begin{IEEEeqnarray}{rCl}
\text{ \textbf{LQR-OPF}:} &&\notag \\
\quad \min_{\mb{S}, \, \mb {z}^s=\{\mb x^s, \mb a^s, \mb u^s\}} \: && \;\;\;\; c(\mb {z}^s)+\frac{T_{\mr{lqr}}}{2}  \gamma \label{eqn:lqropfobj} \IEEEeqnarraynumspace\\
\text{subj.~to} \ && \eqref{eqn:linearssg}, \: \eqref{eqn:linearssh} , \:\mb{a}^s \in \mc{A} \IEEEeqnarraynumspace  \label{eqn:lqropfextra}\\
\IEEEeqnarraymulticol{3}{c}{ \bmat{-\gamma  & (\mb {x}^s-\mb {x}^0)^{\top} \\ \mb {x}^s - \mb {x}^0 & -\mb{S}} \preceq \mb {O}} \IEEEeqnarraynumspace\label{eqn:lqropflmi1}\\
\IEEEeqnarraymulticol{3}{c}{\bmat{ \begin{IEEEeqnarraybox}[][c]{rCl}  \mb{A}^{\top}\mb{S} &+& \mb{S} \mb{A}  \\
			+ \: \mb{B} \mb{Y}&+& \mb{Y}^{\top} \mb{B}^{\top} \IEEEstrut \end{IEEEeqnarraybox} & \mb{S} & \mb{Y}^{\top} \\ \mb{S} & -\mb{Q}^{-1}(\mb{z}^s) & \mb{O} \\
	\mb{Y} & \mb{O} & -\mb{R}^{-1}(\mb{z}^s)} \preceq \mb{O}} \label{eqn:lqropflmi2} \IEEEeqnarraynumspace \\
&&\mb{S} \succeq \mb {O}. \label{eqn:lqropfpd}
\end{IEEEeqnarray}
\end{subequations}
The optimal control law generated from the optimal control problem \eqref{eqngroup:LQR1} is $\Delta \mb {u}'=\mb{K}\Delta \mb{x}'$ where
	{ $\mb{K} = -\mb{R}^{-1}(\mb{z}^{*})\mb{B}(\mb{z}^0)^{\top} \mb{S}^{{*}^{-1}}$} is the optimal state feedback control gain and $(\mb{S}^*,\mb{z}^{*})$ is  an optimizer of \eqref{eqngroup:lqropf}. The reader is referred to~\cite{boydLMI} for the derivation of this control law.

The LQR-OPF problem \eqref{eqngroup:lqropf} is an SDP which solves jointly for the new steady-state variables $\mb{z}^s$ and the matrix $\mb{S}$ while guaranteeing the stability of the linearized system in \eqref{eqn:newlineardynamical}.  The two problems of OPF and stability control are coupled in~\eqref{eqngroup:lqropf} through the LQR weight matrices---with $\mb{Q}^{-1}$ and $\mb{R}^{-1}$ being affine functions---as well as through the term  $\mb{x}^s-\mb{x}^0$ in LMI~\eqref{eqn:lqropflmi1}. It is worth noting that $\gamma$ in~\eqref{eqngroup:lqropf} is equal to the integral in the objective of~\eqref{eqngroup:LQR1} with $t_f=\infty$.

The optimal steady-state solution $\mb{z}^{s}$ of the LQR-OPF problem \eqref{eqngroup:lqropf} satisfies the linearized steady-state equations \eqref{eqngroup:linearss}. In order to obtain an equilibrium for the true nonlinear DAEs \eqref{eqngroup:loadflow1},  generator setpoints from $\mb{z}^s$ are extracted as follows. For generator $i \in \mc{G}$,  set $v_i^{\mr{eq}}=v_{i}^s$. If $i$ is nonslack,  set  $p_{g_i}^{\mr{eq}}=p_{g_i}^s$. If $i$ is slack,  set $\theta_i^{\mr{eq}}=\theta_i^s$.  Then, the system of equations~\eqref{eqngroup:loadflow1} is solved to obtain an equilibrium $\mb{z}^{\mr{eq}}$ for the nonlinear DAEs in~\eqref{eqngroup:DAEs}. This task essentially amounts to solving two separate sets of nonlinear equalities: first  a standard load-flow is solved to obtain the remainder of algebraic variables $q_{g_i}^{\mr{eq}}$ for $i \in \mc{G}$, $v_i^{\mr{eq}}$ for $i \in \mc{L}$, $p_{g_i}^{\mr{eq}}$ if $i$ is slack, and $\theta_{i}^{\mr{eq}}$ for $i$  nonslack. Second, by incorporating  the algebraic variables in generator equations, the equilibrium states and controls of the generators, that is $\mb{x}^{\mr{eq}}$ and $\mb{u}^{\mr{eq}}$, are obtained. This process is detailed in Section~\ref{sec:numtests} that includes simulations.

\begin{figure}[t]
	\centering 
	\scalebox{1.1}{
		\hspace{-0.35cm}\begin{tikzpicture}[scale=1, auto, >=stealth']
		\scriptsize
		
		\matrix[ampersand replacement=\&, row sep=0.8cm, column sep=0.8cm] {
			\& \node[block] (LQR-OPF) {LQR-OPF \eqref{eqngroup:lqropf}}; \& 
			\node[block] (LQR) {LQR};  \&   \node[block] (Observer) {Observer};  \\ 
			\node[block](forecast) {Forecast}; 
			\&  \node[block][align=center] (linear) {Linearized \\ System  \eqref{eqngroup:genericlinear}}; 
			\&   \node[block] [align=center] (nonlinear) {Nonlinear \\ System \eqref{eqngroup:DAEs}}; \&
			\\
		};
		
		\draw [connector] (LQR-OPF) -- node[pos=0.5] {$ \mathbf{x}^{\mr{eq}},\mathbf{u}^{\mr{eq}}$} (LQR);
		\draw [connector] (LQR) -- node {$\mathbf{K}(\hat{\mb{x}}-\mb{x}^{\mr{eq}})$} (nonlinear);
		\draw [rightconnector] (LQR) -- node {$\hat{\mathbf{x}}$} (Observer);
		\draw [connector] (linear) -- (LQR-OPF);
		\draw  ($(linear.north)+(0.3cm,0.45cm)$)  node {$\mb{z}^0$}(LQR-OPF);
		\draw  ($(linear.north)-(0.5cm,-0.60cm)$)  node {$\mb{A}(\mb{z}^0)$}(LQR-OPF);
		\draw   ($(linear.north)-(0.5cm,-0.20cm)$)  node {$\mb{B}(\mb{z}^0)$}(LQR-OPF);
		\draw [connector] (nonlinear) -| node[pos=0.2] {$\mathbf{y}$} (Observer);
		\draw [rightconnector] (linear) --node[pos=0.5] {$\mb{z}^0$} (nonlinear);
		\draw [connector] (forecast.north) |- node[pos=0.5] {$\mathbf{d}^s$} (LQR-OPF.west);  
		\end{tikzpicture}}
	\caption{Overall design of the proposed method. The LQR-OPF uses the forecast information to determine the optimal setpoints for the system. The feedback control law generated by the LQR drives the system to the optimal steady-state with minimized control costs.}
	\label{fig:overalldesign}
\end{figure}
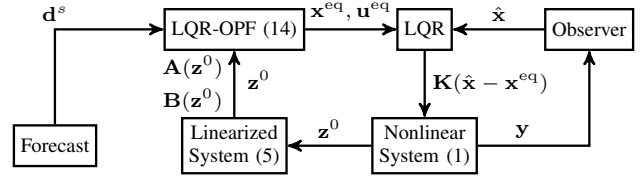  

Figure ~\ref{fig:overalldesign} demonstrates the overall design and the integration of the proposed LQR-OPF into the power system operation. Initially, the system is operating at a known equilibrium point $\mb{z}^0$.    Using this known operating point, the state-space matrices $\mb{A}(\mb{z}^0)$ and $\mb{B}(\mb{z}^0)$ are computed according to \eqref{eqn:newlineardynamical}.  This information is then provided to the LQR-OPF block together with a forecast $\mb{d}^s$ of the demand at the next steady state. The LQR-OPF calculates the optimal generator control setpoints included in $\mb{u}^{\mr{eq}}$ as well as the optimal state setpoints for the next steady state $\mb{x}^{\mr{eq}}$. Upon providing the LQR block with  optimal setpoints $\mb{x}^{\mr{eq}}$ and $\mb{u}^{\mr{eq}}$, the LQR control law $\mb{K}(\mb{x}-\mb{x}^{\mr{eq}})+\mb{u}^{\mr{eq}}$ is then applied to the nonlinear power system and drives the system to the new desired optimal point~$\mb{z}^{\mr{eq}}$.

 The purpose of LQR-OPF~\eqref{eqngroup:lqropf} is to generate \textit{steady-state operating points} with desirable stability properties and minimal secondary control effort, so that the optimal $\mb z^{\mr{eq}}$ in~\eqref{eqngroup:lqropf} is cognizant of dynamic stability constraints.  The proposed LQR-OPF facilitates the development of a more general framework where various dynamical system applications can be integrated with static operational constraints, leading to the computation of stability-aware operating points---more concrete pointers to this end are provided in Section~\ref{sec:conclusion}.

\begin{remark}
In the previous sections, availability of the full state-vector $\mb x(t)$ (i.e., angles, frequencies, and EMFs for generators) is assumed for the optimal control law $\mb u(t)$. This assumption entails that sensors such as phasor measurement units (PMU) generate $\mb x(t)$ for all generator buses in the network. However, a PMU is not needed on all generator buses, as dynamic state estimation tools for power networks can be leveraged. See e.g.,~\cite{cui2015particle,valverde2011unscented} for estimation through Kalman filters and \cite{taha2016risk} for deterministic observers.  For example, it has been previously demonstrated  in~\cite{taha2016risk} that 12 PMUs are sufficient to estimate the 10th order dynamics of the 16 generators in the New England 68-bus power network. 
\end{remark}

\begin{remark}
Problem \eqref{eqngroup:lqropf} can be straightforwardly extended to include convex relaxations of the nonlinear power flow equations (such as SDP relaxations \cite{Madani2015}) in place of the linearized power flow constraints \eqref{eqn:linearssg} and \eqref{eqn:linearssh}. 
\end{remark}

\section{Approximate Solver for the LQR-OPF}
The LQR-OPF problem presented in \eqref{eqngroup:lqropf} has two LMIs and may be computationally expensive to solve for very large networks via interior-point methods.  In this section,  an algorithm to approximately solve \eqref{eqngroup:lqropf} is developed.  

To this end, consider the following nonconvex optimization: 
\begin{subequations}
	\label{eqngroup:ALQROPFIntro}
			\begin{IEEEeqnarray}{rCl}
		\min_{\substack{ \mb{z}^s=\{\mb{x}^s, \mb{a}^s, \mb{u}^s\} \\ \mb{P}}} \: &&  c(\mb{a}^s)+ \frac{T_{\mr{lqr}}}{2} (\mb{x}^s-\mb{x}^0)^\top \mb{P} (\mb{x}^s- \mb{x}^0) 
		\label{eqn:ALQROPFIntroObj} \IEEEeqnarraynumspace \\
		\text{subj.~to} \ && \eqref{eqn:linearssg}, \: \eqref{eqn:linearssh} , \:\mb{a}^s \in \mc{A} \IEEEeqnarraynumspace  \label{eqn:ALQROPFIntroExtra}\\
		\IEEEeqnarraymulticol{3}{c}{\mb{A}^{\top} \mb{P}+ \mb{P}\mb{A}-\mb{P}\mb{B}\mb{R}^{-1}(\mb{z}^s)\mb{B}\mb{P}+\mb{Q}(\mb{z}^s)=\mb{O},  \label{eqn:ALQROPFIntroCare}}
		\end{IEEEeqnarray}
\end{subequations}
where \eqref{eqn:ALQROPFIntroCare} is the well-known continuous algebraic Riccati equation (CARE). Problem~\eqref{eqngroup:ALQROPFIntro} acts as a surrogate for~\eqref{eqngroup:lqropf} because the optimal values of problems~\eqref{eqngroup:lqropf} and~\eqref{eqngroup:ALQROPFIntro} are equal. To see this, note first that the objective $c(\mb{a}^s)$ and constraints  \eqref{eqn:linearssg}, \eqref{eqn:linearssh}, and $\mb{a}^s \in \mc{A}$ are the same for both problems. Then, using the theory in \cite[Sec.~V]{balakrishnan2003semidefinite}, it can be shown that for any solution $\mb{P}$ of the CARE in~\eqref{eqn:ALQROPFIntroCare}, the quadratic form $(\mb{x}^s-\mb{x}^0)^\top \mb{P} (\mb{x}^s- \mb{x}^0)$ is equal to $\gamma$ in~\eqref{eqn:lqropfobj}.

Although~\eqref{eqngroup:ALQROPFIntro} is nonconvex, an efficient algorithm that alternates between variables $\mb{z}^s$ and $\mb{P}$ can be used to approximately solve~\eqref{eqngroup:ALQROPFIntro}. Specifically, with $\mb{z}^s$ fixed,  the CARE constraint \eqref{eqn:ALQROPFIntroCare} can be solved for  to obtain $\mb{P}$. Then,  plugging  $\mb{P}$ into  \eqref{eqngroup:ALQROPFIntro} and removing \eqref{eqn:ALQROPFIntroCare}  yields a quadratic program (QP), which can be efficiently solved. The process can be repeated as long as the objective in~\eqref{eqn:ALQROPFIntroObj} is improved, and is summarized in Algorithm~\ref{algorithm:ALQR-OPF}. The algorithm is referred to as ALQR-OPF. 
 
The technology for solving the CARE and QPs with very large number of variables has significantly matured and enables the efficient implementation of Algorithm~\ref{algorithm:ALQR-OPF}. The numerical tests of Section \ref{sec:numtests} indicate that Algorithm~\ref{algorithm:ALQR-OPF} does not practically compromise optimality and is scalable to networks with thousands of buses. 

\begin{algorithm}[t]
	\caption{Approximate solver for LQR-OPF (ALQR-OPF)}
	\label{algorithm:ALQR-OPF}
	\begin{algorithmic}[1]
	    \State Initialize $o^{\mr{best}} \gets +\infty$
	    \State Set $\mb{P}^{(0)}$ to the solution of CARE \eqref{eqn:ALQROPFIntroCare} for $\mb{z}^s=\mb{z}^0$
	    \For{$k=1:k^{\max}$}
	    \State Set $\mb{z}^{(k)}=\{\mb{x}^{(k)}, \mb{a}^{(k)}, \mb{u}^{(k)}\}$ to the optimal solution of \eqref{eqn:ALQROPFIntroObj}--\eqref{eqn:ALQROPFIntroExtra} for $\mb{P} = \mb{P}^{(k-1)}$
	    \State Set $\mb{P}^{(k)}$ to the solution of CARE \eqref{eqn:ALQROPFIntroCare} for $\mb{z}^s= \mb{z}^{(k)}$
	    \If{  $c(\mb{a}^{(k)})+ \frac{T_{\mr{lqr}}}{2} (\mb{x}^{(k)}-\mb{x}^0)^\top \mb{P}^{(k)} (\mb{x}^{(k)}- \mb{x}^0)  ) < o^{\mr{best}}$ }
  			\State $ o^{\mr{best}} \gets c(\mb{a}^{(k)})+ \frac{T_{\mr{lqr}}}{2} (\mb{x}^{(k)}-\mb{x}^0)^\top \mb{P}^{(k)} (\mb{x}^{(k)}- \mb{x}^0)$ 
  			\State $\mb{z}^s=\{\mb{x}^s, \mb{a}^s, \mb{u}^s\} \gets \mb{z}^{(k)}=\{\mb{x}^{(k)}, \mb{a}^{(k)}, \mb{u}^{(k)}\}$
  		    \State $ \mb{P}^*= \mb{P}^{(k)}$
    \EndIf
	    \EndFor
	\end{algorithmic}
\end{algorithm}

\color{black}
\section{Generator Model, Power Flows, and OPF}
\label{sec:genmodelpfopf}
The developments in this paper and the proposed LQR formulation are applicable to any generator dynamical model with control inputs.  However,  as an example, a  specific form of mappings $\mb{g}$ and $\mb{h}$ is given in this section, which will be used for the numerical tests of Section~\ref{sec:numtests}.

\subsection{Generator model}
\label{subsec:genmodel}
The fourth order model of the synchronous generator internal dynamics for node $i \in \mathcal{G}$ can be written as
\vspace{-0.2cm}
 {\begin{subequations}
 \label{eqngroup:gendiff}
\begin{IEEEeqnarray}{rCl}
\dot{\delta}_i&=& \omega_i - \omega_s \label{eqn:deltadynamics} \\
 \dot{\omega}_i &=& \frac{1}{M_i}\left[m_i - D_i ( \omega_i - \omega_s) - p_{g_i}\right] \label{eqn:omegadynamic} \\
\dot{e}_i &=& \frac{1}{\tau_{d_i}}\left[-\frac{x_{di}}{x_{di}'} e_i + \frac{x_{di}-x_{di}'}{x_{di}'} v_i \cos(\delta_i- \theta_i) + f_i\right] \IEEEeqnarraynumspace \label{eqn:edynamics}  \\
\dot{m}_i&=& \frac{1}{\tau_{c_i}} \left[ r_i(t)- \frac{1}{R_i} (\omega_i - \omega_s)- m(t)\right] \label{eqn:mdynamics}
\end{IEEEeqnarray}
\end{subequations}}where $M_i$ is the rotor's inertia constant ($\mr{pu} \times \mr{sec}^2$), $D_i$ is the damping coefficient ($\mr{pu} \times \mr{sec}$), $\tau_{d_i}$ is the direct-axis open-circuit time constant ($\mr{sec}$), $x_{d_i}$ and $x_{q_i}$ are respectively the direct- and quadrature-axis synchronous reactances, and $x'_{di}$ is the direct-axis transient reactance ($\mr{pu}$).   Equation \eqref{eqn:mdynamics} is a simplified model of a prime-mover generator with $\tau_{c_i}$ as the charging time ($\mr{sec}$) and  a speed-governing  mechanism with  $R_i$ as the regulation constant ($\frac{ \mr{Rad} \times \mr{Hz}}{\mr{pu}}$).
 The mapping $\mb{g}$ defined in \eqref{eqn:genericdifferential} is given by concatenating  \eqref{eqngroup:gendiff} for $i \in \mc{G}$.

The following algebraic equations relate the generator real and reactive power output with generator voltage, internal EMF, and internal angle and must hold at any time instant for generator nodes $i \in \mc{G}$: 
{
\begin{subequations}
\label{eqngroup:genalg}
\begin{IEEEeqnarray}{rCl}
0 &=&-p_{g_i}+\frac{e_iv_i}{x_{di}'} \sin(\delta_i - \theta_i) \notag\\
&&+\frac{x_{di}'-x_{q_i}}{2x_{qi}x_{di}'} v_i^2 \sin[2(\delta_i - \theta_i)] \IEEEeqnarraynumspace \label{eqn:genalgpg} \\
0&=& -q_{g_i}+ \frac{e_iv_i}{x_{di}'} \cos (\delta_i - \theta_i) -\frac{ x'_{di}+ x_{qi}}{2x_{qi}x'_{di}}v_i^2 \notag\\ &&+ \frac{  x'_{di}- x_{qi}}{2x_{qi}x'_{di}}v_i^2\cos[2(\delta_i - \theta_i)]. \IEEEeqnarraynumspace \label{eqn:genalgqg}
\end{IEEEeqnarray}
\end{subequations}
}

\subsection{Power flow equations}
\label{subsec:powerflowequations}
Let $\mb{Y}=\boldsymbol{\mathscr{G}}+j\boldsymbol{\mathscr{B}}$ denote the network bus admittance matrix based on the $\pi$-model of transmission lines.  Notice that $\mb{Y}$  may include transformers, tap-changing voltage regulators, and phase shifters \cite{smartgridtutorial2013}. The power flow equations are
{ \begin{subequations}
\label{eqngroup:powerflows}
\begin{IEEEeqnarray}{rCl}
-p_{l_i} &=&- p_{g_i}+ \: \mathscr{G}_{i,i} v_i^2  +  \sum\nolimits_{j \in \mc{N}_i} \left[ \mathscr{G}_{i,j} v_i v_j \cos\theta_{ij} \right.  \notag \\
&&  \mspace{110mu} + \:  \left. \mathscr{B}_{i,j} v_i v_j \sin\theta_{ij} \right], i \in \mc{G},\label{eqn:pgpf} \\
-q_{l_i} &=&-q_{g_i} - \: \mathscr{B}_{i,i} v_i^2  +   \sum\nolimits_{j \in \mc{N}_i} \left[\mathscr{G}_{i,j} v_i v_j \sin\theta_{ij} \right. \notag \\
&&  \mspace{110mu} - \: \left.  \mathscr{B}_{i,j} v_i v_j \cos\theta_{ij}\right], i \in \mc{G},  \label{eqn:qgpf}\\
-p_{l_i} &=&  \: \mathscr{G}_{i,i} v_i^2  +  \sum\nolimits_{j \in \mc{N}_i} \left[ \mathscr{G}_{i,j} v_i v_j \cos\theta_{ij} \right.  \notag \\
&&   \mspace{110mu} + \:  \left. \mathscr{B}_{i,j} v_i v_j \sin\theta_{ij} \right], i \in \mc{L},\label{eqn:plpf} \\
-q_{l_i} &=& - \: \mathscr{B}_{i,i} v_i^2  +   \sum\nolimits_{j \in \mc{N}_i} \left[\mathscr{G}_{i,j} v_i v_j \sin\theta_{ij} \right. \notag \\
&&  \mspace{110mu} - \: \left.  \mathscr{B}_{i,j} v_i v_j \cos\theta_{ij}\right], i \in \mc{L},  \label{eqn:qlpf}
\end{IEEEeqnarray}
\end{subequations}}where $\theta_{ij}:=\theta_i - \theta_j$; $p_{l_i}:=p_{l_i}(t)$  and  $q_{l_i}:=q_{l_i}(t)$ are respectively the real and reactive power demands at node $i$ modeled as a time-varying constant-power load, i.e.,  $p_{l_i}$ and $q_{l_i}$ are not functions of $v_i$. 
The mapping $\mb{h}$ in \eqref{eqn:genericalgebraic} is given by concatenating \eqref{eqngroup:genalg} for $i \in \mc{G}$ and \eqref{eqngroup:powerflows} for $i \in \mc{N}$. By defining $\mb{p}_{l_{\mc{G}}}=\{p_{l_i}\}_{i \in \mc{G}}$, $\mb{q}_{l_{\mc{G}}}=\{q_{l_i}\}_{i \in \mc{G}}$, $\mb{p}_{l_{\mc{L}}}=\{p_{l_i}\}_{i \in \mc{L}}$, $\mb{q}_{l_{\mc{L}}}=\{q_{l_i}\}_{i \in \mc{L}}$,  the vector $\mb{d}= \{\mb{0}_{2G}, -\mb{p}_{l_\mc{G}}, -\mb{q}_{l_\mc{G}}, -\mb{p}_{l_\mc{L}}, -\mb{q}_{l_\mc{L}}\}$ is obtained. 

\subsection{Optimal power flow}
\label{subsec:opf}
The standard OPF problem typically only considers the algebraic variable $\mb{a}^s$ and is given as 
\begin{IEEEeqnarray}{lrCl}
\min_{\mb{a}^s} c( \mb{a}^s) \;\;\;\: 
\text{subj. to} \;\;\;\: \eqref{eqngroup:powerflows}
\: \text{ and } \mb{a}^s \in \mc{A}. \label{eqn:opf}
\end{IEEEeqnarray}
 If the cost function $c(\mb{z}^s)$ of the generalized OPF  in \eqref{eqngroup:generalopf} is only a function of algebraic variables $\mb{a}^s$,  \eqref{eqngroup:generalopf} can be solved by solving the standard OPF problem  \eqref{eqn:opf} to obtain the optimal $\mb{a}^s$. The variables $\mb{x}^s$ and $\mb{u}^s$ can then be found by plugging in $\mb{a}^s$ in equations \eqref{eqngroup:genalg} and \eqref{eqngroup:gendiff} while setting $\dot{\mb{x}}=\mb{0}$.
 
 \section{Numerical Simulations}
\label{sec:numtests}
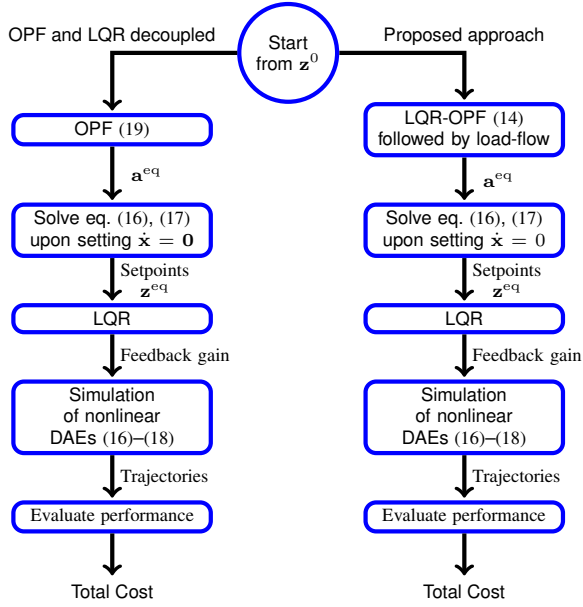
\begin{figure}[t]\scriptsize
	\centering
	\begin{tikzpicture} [
	auto,
	block/.style    = { rectangle, draw=blue, ultra thick, text width=10em, text centered,
		rounded corners, minimum height=1em },
	block1/.style    = { circle, draw=blue, ultra thick, text width=4em, text centered,
		rounded corners, minimum height=1em },
	line/.style     = { draw, ultra thick, ->, shorten >=1pt },
	]
	
	\matrix [column sep=20mm, row sep=6mm] {
		\node [block] (SAOPF2) {\textsf{OPF \eqref{eqn:opf} } } ;  & \node [block] (SAOPF1) {\textsf{LQR-OPF \eqref{eqngroup:lqropf} followed by load-flow}};    \\
		\node [block] (GE2) {\textsf{Solve eq.~\eqref{eqngroup:gendiff}, \eqref{eqngroup:genalg} upon setting $\dot{\mb{x}}=\mb{0}$}}; 	& \node[block] (GE1){\textsf{Solve eq. \eqref{eqngroup:gendiff}, \eqref{eqngroup:genalg} upon setting $\dot{\mb{x}}=0$}}; \\
		\node [block] (SolveLQR2) {\textsf{LQR}}; 	& \node [block] (SolveLQR1) {\textsf{LQR}};             \\
		\node [block] (DAE2) {\textsf{Simulation of nonlinear DAEs \eqref{eqngroup:gendiff}--\eqref{eqngroup:powerflows}
		}}; 	& \node [block] (DAE1) {\textsf{Simulation of nonlinear DAEs \eqref{eqngroup:gendiff}--\eqref{eqngroup:powerflows}}};  \\
		\node [block] (ObtainP2){Evaluate  performance}; & \node [block] (ObtainP1){Evaluate  performance};         \\
		\node [text centered] (TC2) {\textsf{Total Cost}};	& \node [text centered] (TC1) {\textsf{Total Cost}};         \\
	};
	\begin{scope}[xshift=0cm,yshift=4cm]	
	\matrix [column sep=5mm, row sep=3mm] {
		& \node [block1] (SFS) {\textsf{Start from} $\mb{z}^0$};   & \\
	};\end{scope}
	\begin{scope} [every path/.style=line]
	\path (SFS) -| node [above] {\textsf{Proposed approach}} (SAOPF1);
	\path (SFS) -| node [above] {\textsf{OPF and LQR decoupled}} (SAOPF2);
	\path (SAOPF1)      --    node [midway] {\: $\mb{a}^{\mr{eq}}$} (GE1);
	\path (SAOPF2)      --    node [midway] {\: $\mb{a}^{\mr{eq}}$} (GE2);
	\path (GE1)   --    node [align=center, midway] {Setpoints \\ $\mb{z}^{\mr{eq}}$} (SolveLQR1);
	\path (GE2)   --    node [align=center, midway] {Setpoints \\ $\mb{z}^{\mr{eq}}$} (SolveLQR2);
	\path (SolveLQR1)    --    node [midway] {Feedback gain} (DAE1);
	\path(DAE1) -- node [midway] {Trajectories} (ObtainP1);
	\path (ObtainP1) --    node [midway] {} (TC1);
	\path (SolveLQR2)    --    node [midway] {Feedback gain} (DAE2);
	\path(DAE2) -- node [midway] {Trajectories} (ObtainP2);
	\path (ObtainP2) --  node [midway] {}  (TC2);
	\end{scope}
	\end{tikzpicture}
	\caption{Diagram showing the simulation steps. }
	\label{fig:Diagram}
	\vspace{-0.3cm}	
\end{figure}

This section provides a numerical assessment of the advantages of the  LQR-OPF  in comparison to a method where OPF and load-following control problems are treated separately. Prior to analyzing the case studies, the general simulation workflow of Fig.~\ref{fig:Diagram} is described. 

The decoupled approach, one where OPF and stability are solved separately, is considered on the left-hand side of Fig.~\ref{fig:Diagram}.
Initially, the system is in steady-state and operates at $\mb{z}^0$.  For a forecasted load demand, $\mb{d}^s = \mb{d}^0 + \Delta \mb{d}^s$,  the OPF \eqref{eqn:opf} is solved yielding  optimal steady-state setpoints $\mb{a}^s$,  including generator real and reactive powers ($p_{g_i}^s, q_{g_i}^s$) and  generator  voltage magnitudes and phases ($v_i^s, \theta_i^s$). For OPF, it holds that $\mb{a}^{\mr{eq}}=\mb{a}^s$.  The algebraic variables are then utilized  to solve  \eqref{eqngroup:gendiff}  upon setting $\dot{\mb{x}}=\mb{0}$ together with~\eqref{eqngroup:genalg} to obtain steady-state setpoints of generator states $\mb{x}^\mr{eq}$ and control inputs $\mb{u}^\mr{eq}$.  The next steady-state equilibrium is then simply represented as $\mb{z}^{\mr{eq}}=(\mb{x}^{\mr{eq}}, \mb{a}^{\mr{eq}}, \mb{u}^\mr{eq})$. The DAEs~\eqref{eqn:genericdifferential}--\eqref{eqn:genericalgebraic}, upon being subjected to the new load $\mb{d}^s$,   depart from the initial equilibrium $\mb{z}^0$.  To steer the DAEs \eqref{eqn:genericdifferential}--\eqref{eqn:genericalgebraic} to the next desired equilibrium point $\mb{z}^\mr{eq}$, LQR is performed. Dynamic performance as well as costs of  steady-state and load-following control are evaluated. The  standard OPF \eqref{eqn:opf}  is solved using MATPOWER's \texttt{runopf.m}.

The proposed  methodology is considered on the right-hand side of  Fig.~\ref{fig:Diagram} where  the OPF block is replaced by  LQR-OPF \eqref{eqngroup:lqropf} followed by a load-flow.  LQR-OPF  obtains  optimal generator voltage setpoints $\{v_{i}^s\}_{i \in \mc{G}}$, real power setpoints $\{p_{g_i}^s\}_{i \in \mc{G} \setminus \{i_\mr{slack}\}}$,  and $\theta_{i_{\mr{slack}}}^s$, while accounting for  load-following costs that drive the DAEs \eqref{eqn:genericdifferential}--\eqref{eqn:genericalgebraic}  to the next desired equilibrium.  These obtained setpoints are then input to a standard load-flow (performed by MATPOWER's \texttt{runpf.m}) by setting the following for $i \in \mc{G}$:  $v_i^{\mr{eq}}=v_{i}^s$, $p_{g_i}^{\mr{eq}}=p_{g_i}^s$ for nonslack $i$, and   $\theta_{i_{\mr{slack}}}^{\mr{eq}}=\theta_{i_{\mr{slack}}}^s$. This process yields the remaining algebraic variables  in $\mb{a}^\mr{eq}$. Similar to the previous approach, the DAEs~\eqref{eqngroup:gendiff} and \eqref{eqngroup:genalg} are solved after setting $\dot{\mb{x}}=\mb{0}$ yielding optimal steady-state setpoints of states $\mb{x}^{\mr{eq}}$ and controls $\mb{u}^{\mr{eq}}$. The next equilibrium is then given by $\mb{z}^{\mr{eq}}=(\mb{x}^{\mr{eq}}, \mb{a}^{\mr{eq}}, \mb{u}^{\mr{eq}})$. Once the DAE system departs from its initial equilibrium point $\mb{z}^0$ due to $\mb{d}^s$, LQR is applied to drive the system to the  desired equilibrium $\mb{z}^{\mr{eq}}$.  LQR-OPF is solved using the CVX optimization toolbox \cite{cvx1}. 
Finally, a third approach where LQR-OPF in Fig.~\ref{fig:Diagram} is replaced by  ALQR-OPF is considered.  In Algorithm~\ref{algorithm:ALQR-OPF}, the CARE is solved by MATLAB's \texttt{care.m} and the QP using CVX. Comparisons of the three approaches (OPF, LQR-OPF, and ALQR-OPF) are provided in Table~\ref{table:costs}.

The workflow of Fig.~\ref{fig:Diagram} is applied to various networks as described in the next section.  It is worth emphasizing  that even though the derived control law required linearized system dynamics, all simulations are performed on the actual \emph{nonlinear} dynamics and power flow equations, per the DAEs in \eqref{eqn:genericdifferential}--\eqref{eqn:genericalgebraic}. The dynamical system was modeled using {MATLAB's} \texttt{ode} suite. All computations use a personal computer with $32.0~\mr{GB}$  RAM and $3.60~\mr{GHz}$  CPU processor.  MATLAB scripts that simulate the ensuing case studies are provided online \cite{GithubCodes}.
 \subsection{Description of test networks}
The workflow of Fig.~\ref{fig:Diagram} is applied to 9-, 14-, 39-, and 57-bus test networks as well as the 200-bus Illinois network. Steady-state data are obtained from the corresponding  case files available in MATPOWER~\cite{matpower2011}.   These data include bus admittance matrix $\mb{Y}$, steady-state real and reactive power demands, that is, $\mb{p}_l^0:=\{\mb{p}_{l_\mc{G}}^0, \mb{p}_{l_\mc{L}}^0\}$ and $\mb{q}_l^0:=\{\mb{q}_{l_\mc{G}}^0, \mb{q}_{l_\mc{L}}^0\}$,  cost of  real power generation $c(\mb{a}^s)=\sum_i  c_2 (p_{g_i}^s)^2 + c_1 p_{g_i}^s + c_0$, as well as the limits of power generation and voltages in constraint \eqref{eqn:generalopfextraconstraints}.  Machine constants of \eqref{eqn:deltadynamics}--\eqref{eqn:edynamics} and \eqref{eqngroup:genalg} are taken from the PST toolbox, from case files \texttt{d3m9bm.m}, \texttt{d2asbegp.m}, and \texttt{datane.m}, respectively for the 9-, 14-, and 39-bus networks.   For the 57- and the 200-bus networks, as well as the governor model of~\eqref{eqn:mdynamics} for all networks, based on ranges of values provided in PST \cite{PST1992}, typical parameter values of $M_i=0.2~\mr{pu} \times \mr{sec}^2$, $D_i=0~\mr{pu} \times \mr{sec}$, $\tau_{d_i}=5~\mr{sec}$, $x_{d_i}=0.7~\mr{pu}$, $x_{q_i}=0.5~\mr{pu}$, $x'_{d_i}=0.07~\mr{pu}$, $\tau_{c_i}=0.2~\mr{sec}$, and $R_i=0.02~ \frac{\mr{Hz}}{\mr{pu}}$ have been selected.

\begin{table*}[!ht]
	\scriptsize
	\centering
		\renewcommand{\arraystretch}{1.2}
		\caption{Costs Comparison between LQR-OPF, ALQR-OPF, and OPF methods at $\alpha=0.6$ under LQR }
		\label{table:costs}
		\begin{tabular}{|c|c|c|c|c|c|c|c|c|c|c|}
			\hline
			\multirow{3}{*}{Network}&	\multirow{3}{*}{Method}   & \multirow{3}{*}{Obj.} &    Steady-state & Control    &  Total &  Comp. & Control & \multirow{2}{*}{Total}      & Max.&  Max. \\
			& & & cost& est. cost& est. cost & time & cost &  & freq. dev.  &  volt. dev.  \\
			& &  & ($\$$) & ($\$$) &  ($\$$)  & (seconds) & ($\$$) & ($\$$)   & ($\mr{Hz}$) &   ($\mr{pu}$)  \\
			\hline \hline
			\multirow{3}{*}{9-bus}   & LQR-OPF         & 6168.78         & 6144.11         & 26.30           & 6170.41         & 2.57    & 28.02           & 6172.13         & 0.0145          & 0.0262   \\        
			& ALQR-OPF        & 6168.78         & 6144.25         & 26.16           & 6170.41         & 0.85             & 27.87           & 6172.11         & 0.0145          & 0.0261       \\
			& OPF             & --- & 6113.60         & 329.02          & 6442.62         & 0.64   & 316.01          & 6429.61         & 0.0202          & 0.1240      \\
			\hline	
			\multirow{3}{*}{14-bus}   & LQR-OPF         & 9209.27         & 9178.82         & 36.20           & 9215.01         & 2.32       & 36.95           & 9215.77         & 0.0052          & 0.0297     \\      
			& ALQR-OPF        & 9209.28         & 9178.57         & 36.48           & 9215.04         & 0.74            & 37.26           & 9215.83         & 0.0053          & 0.0297   \\
			& OPF             & ---& 9127.35         & 248.29          & 9375.64         & 0.26    & 248.13          & 9375.48         & 0.0120          & 0.0297  \\
			\hline  
			\multirow{3}{*}{39-bus}   & LQR-OPF         & 55560.94        & 52872.28        & 2471.68         & 55343.95        & 18.95           & 3120.73         & 55993.01        & 0.0804          & 0.0820       \\
			& ALQR-OPF        & 55561.20        & 52871.80        & 2472.94         & 55344.74        & 1.39           & 3128.34         & 56000.14        & 0.0804          & 0.0816    \\
			& OPF             & --- & 51386.02        & 12504.70        & 63890.72        & 0.23    & 13468.42        & 64854.44        & 0.1333          & 0.1157    \\
			\hline 
			\multirow{3}{*}{57-bus}       & LQR-OPF         & 50169.04        & 48322.27        & 2468.43         & 50790.70        & 6.01        & 2306.30         & 50628.57        & 0.0602          & 0.0599   \\      
			& ALQR-OPF        & 50177.06        & 48368.33        & 2410.75         & 50779.08        & 0.87             & 2260.24         & 50628.57        & 0.0593          & 0.0600     \\
			& OPF             & --- & 47199.75        & 5829.34         & 53029.09        & 0.22              & 5889.98         & 53089.73        & 0.0944          & 0.0637\\
			\hline
			\multirow{3}{*}{200-bus}  & LQR-OPF         & 52731.74        & 48347.04        & 4090.34         & 52437.38        & 44070.98       & 4089.44         & 52436.48        & 0.0283          & 0.0608  \\
			& ALQR-OPF        & 54575.99        & 50349.65        & 2424.85         & 52774.50        & 2.92  & 2526.84         & 52876.49        & 0.0199          & 0.0589    \\            
			& OPF             & --- & 48271.72        & 8724.08         & 56995.80        & 0.72         & 9258.04         & 57529.76        & 0.0468          & 0.0745        \\
			\hline      
	\end{tabular}
\end{table*}

\subsection{Regulation cost matrices $\mb{Q}$ and $\mb{R}$}
In accordance with~\eqref{eqngroup:generalqrinv}, $\mb{Q}^{-1}$ and $\mb{R}^{-1}$ are selected to be diagonal with affine entries as follows
\begin{subequations}
	\label{eqngroup:QRinv}
	\begin{IEEEeqnarray}{rCl}
	\mb{Q}_{\omega_i, \delta_i, m_i}^{-1} &=& 	\mb{R}_{r_i}^{-1} =\left(1- \alpha \frac{p_{g_i}^s}{p_{g_i}^{\max}}\right),  \label{eqn:Qomega}\\
	\mb{Q}_{e_i}^{-1} &=& 	\mb{R}_{f_i}^{-1}= \left(1- \alpha \frac{q_{g_i}^s}{q_{g_i}^{\max}}\right), \label{eqn:Qe}
	\end{IEEEeqnarray}
\end{subequations}
where $\mb{Q}_{\omega_i,\delta_i, m_i}$ refers to the diagonal entries of $\mb{Q}$ corresponding to $\omega_i$, $\delta_i$, and $m_i$. Matrices $\mb{Q}_{e_i}$,   $\mb{R}_{r_i}$, and $\mb{R}_{f_i}$ are similarly defined.  Parameter $\alpha$ is  in the interval $[0,1)$ that determines the amount of coupling between steady-state quantities and control costs through matrices $\mb{Q}$ and $\mb{R}$. Quantities $p_{g_i}^{\max}$ and $q_{g_i}^{\max}$ are respectively the maximum real and reactive power limits of generator $i \in \mc{G}$. 

The rationale behind choosing \eqref{eqn:Qomega}  is that angle and frequency instability are usually remedied by generating real power.  In this case, increase in steady-state real power generation $p_{g_i}^s$ leads to  higher cost of frequency regulation. Similarly,  the rationale behind choosing \eqref{eqn:Qe} is that voltage stability is typically correlated with reactive power injection. This choice means that an increase in steady-state reactive power generation $q_{g_i}^s$ incurs a higher cost of voltage regulation.

\begin{figure*}[!ht]
	\centering
	\subfloat[]{\includegraphics[scale=0.3]{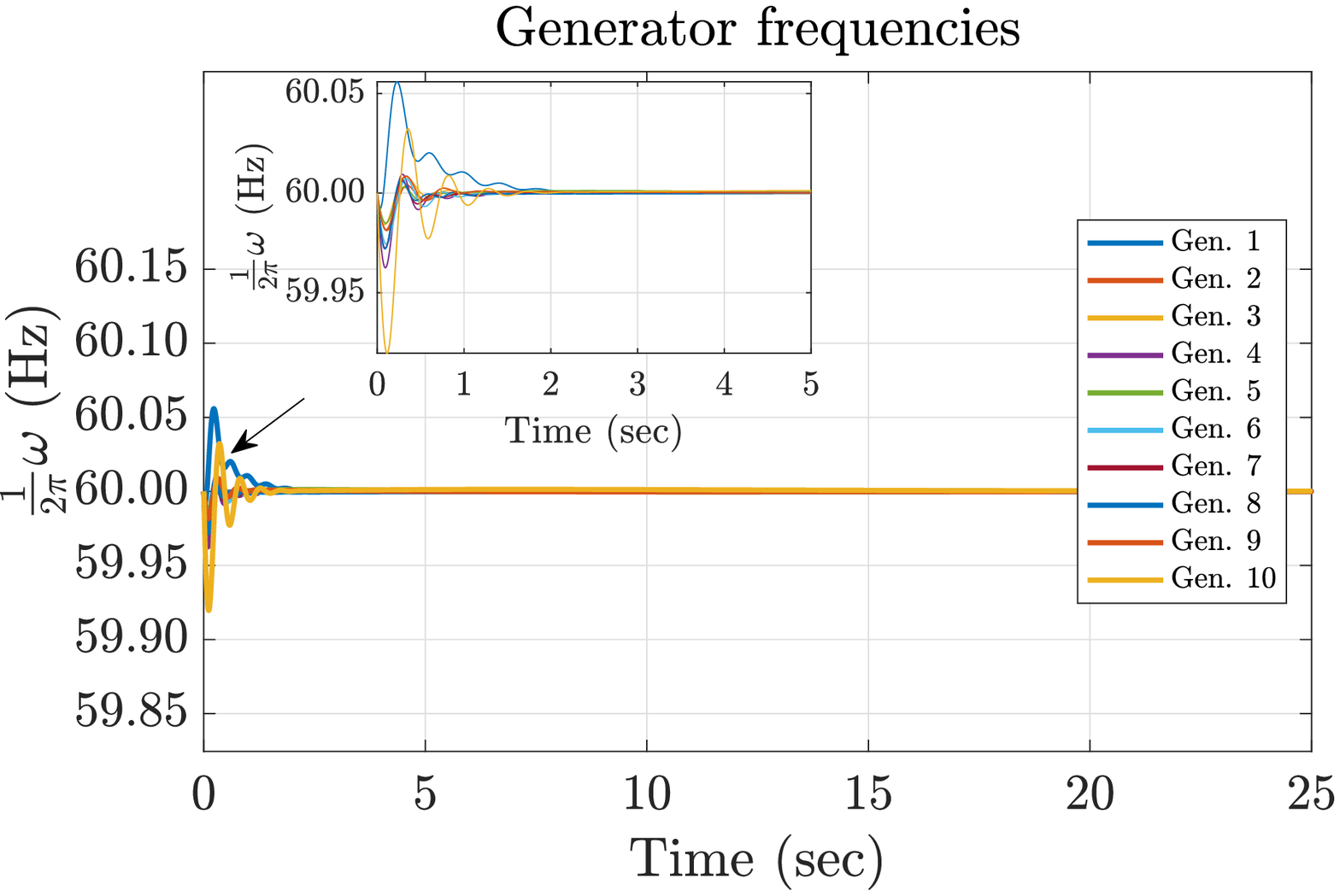}\label{fig:LQR-OPFFreq}}{} \quad
	\subfloat[]{\includegraphics[scale=0.3]{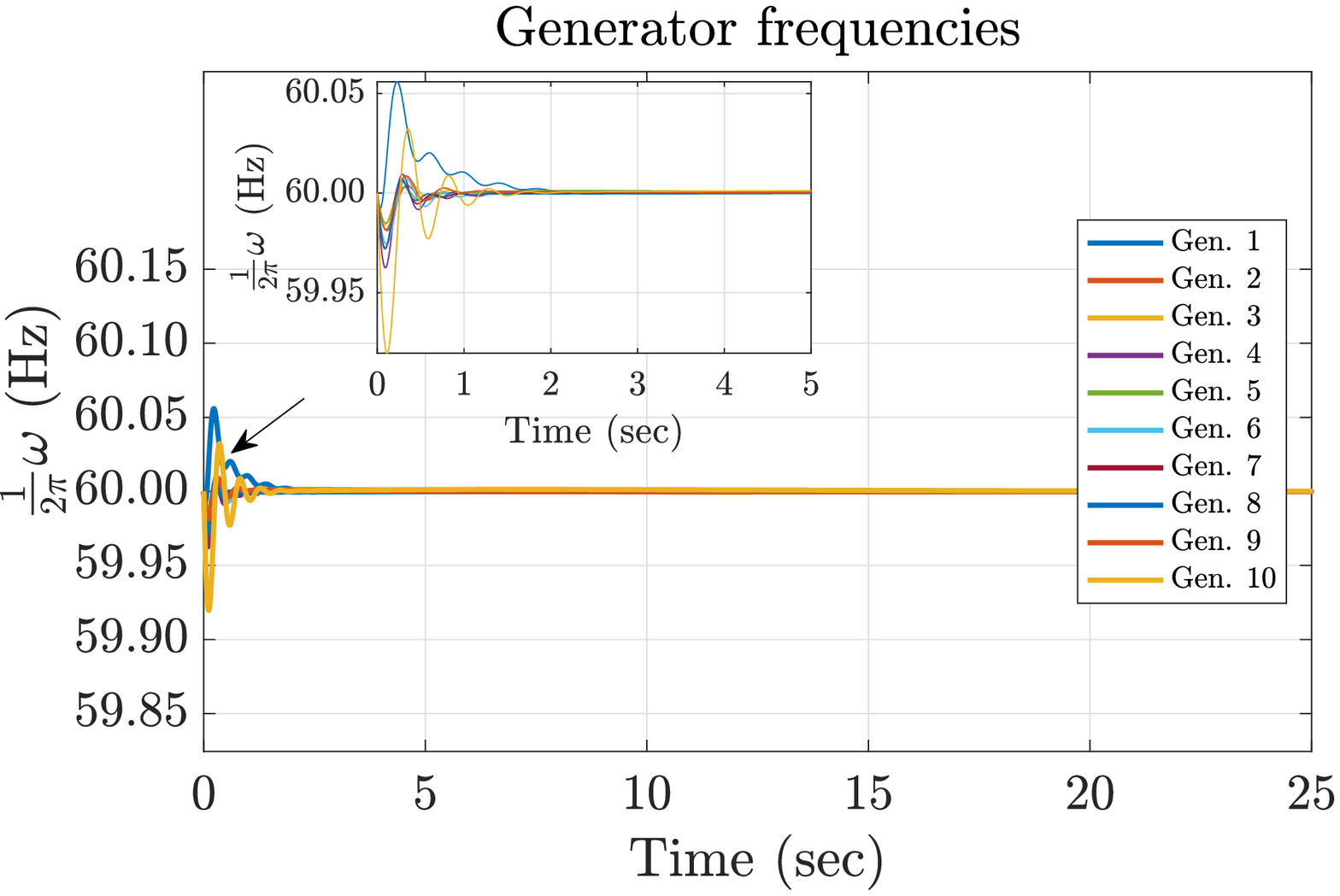}\label{fig:ALQR-OPFFreq}}{}  \quad
	\subfloat[]{\includegraphics[scale=0.3]{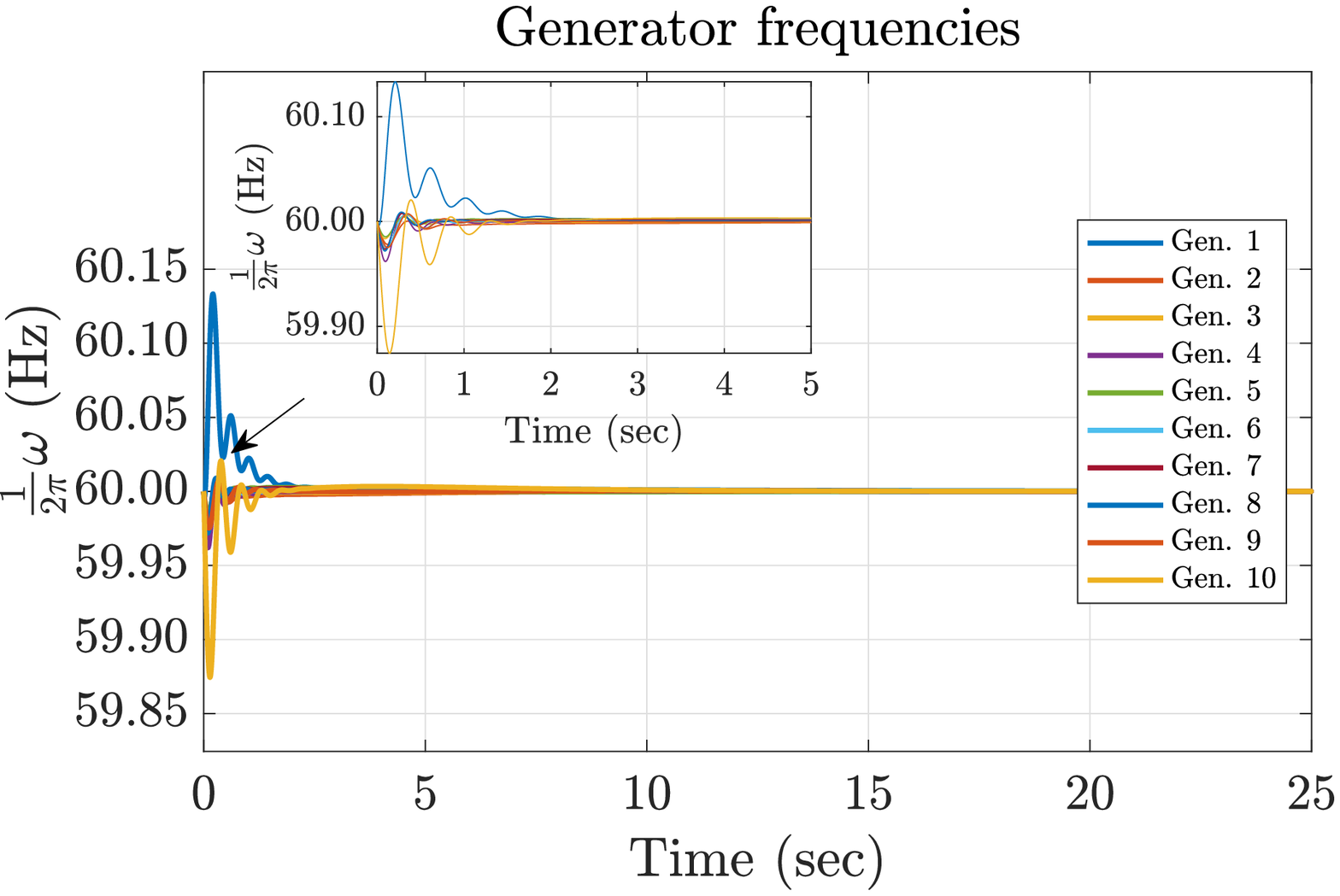}\label{fig:OPFFreq}}{} \\
	\subfloat[]{\includegraphics[scale=0.3]{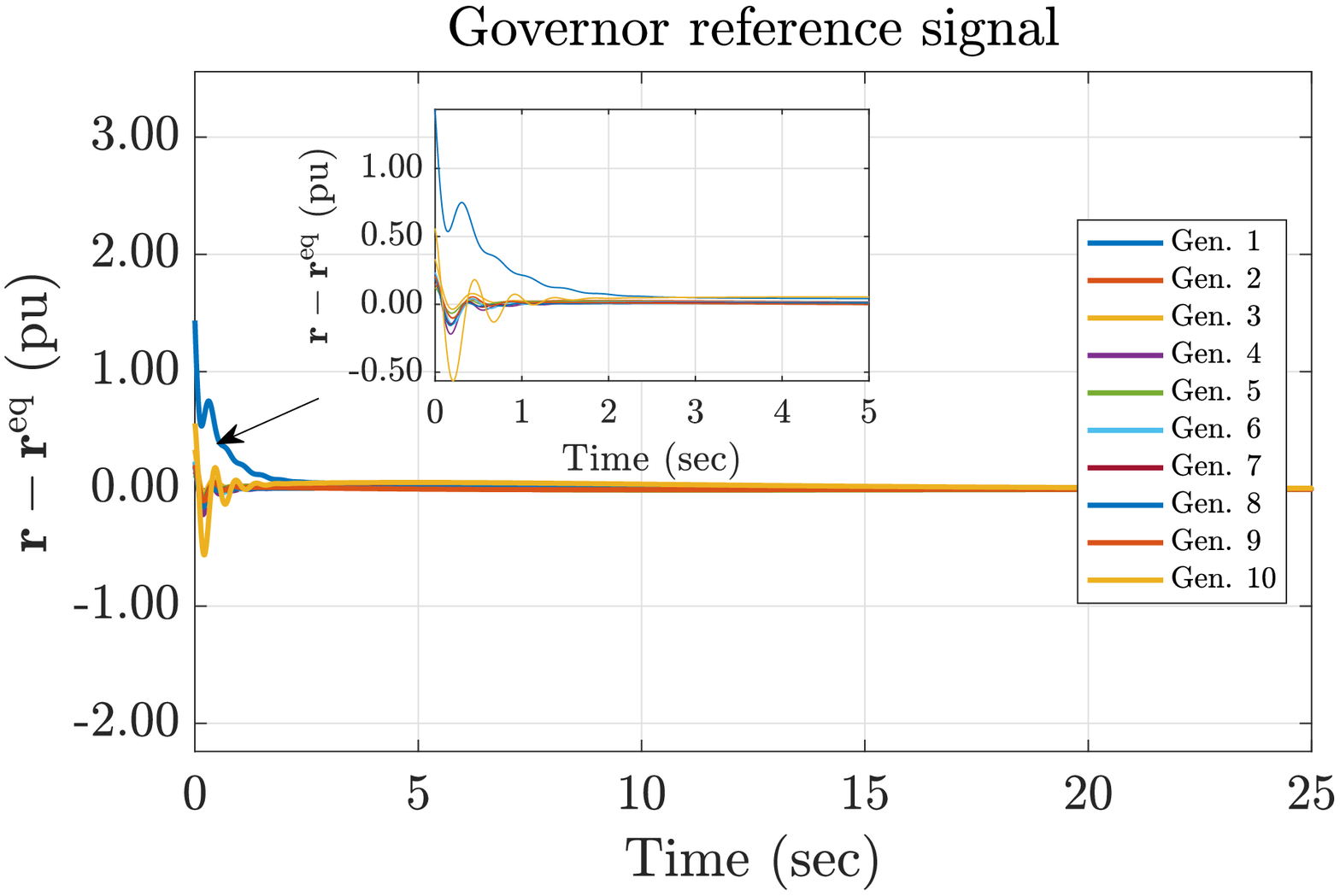}\label{fig:LQR-OPFr}}{}  \quad
	\subfloat[]{\includegraphics[scale=0.3]{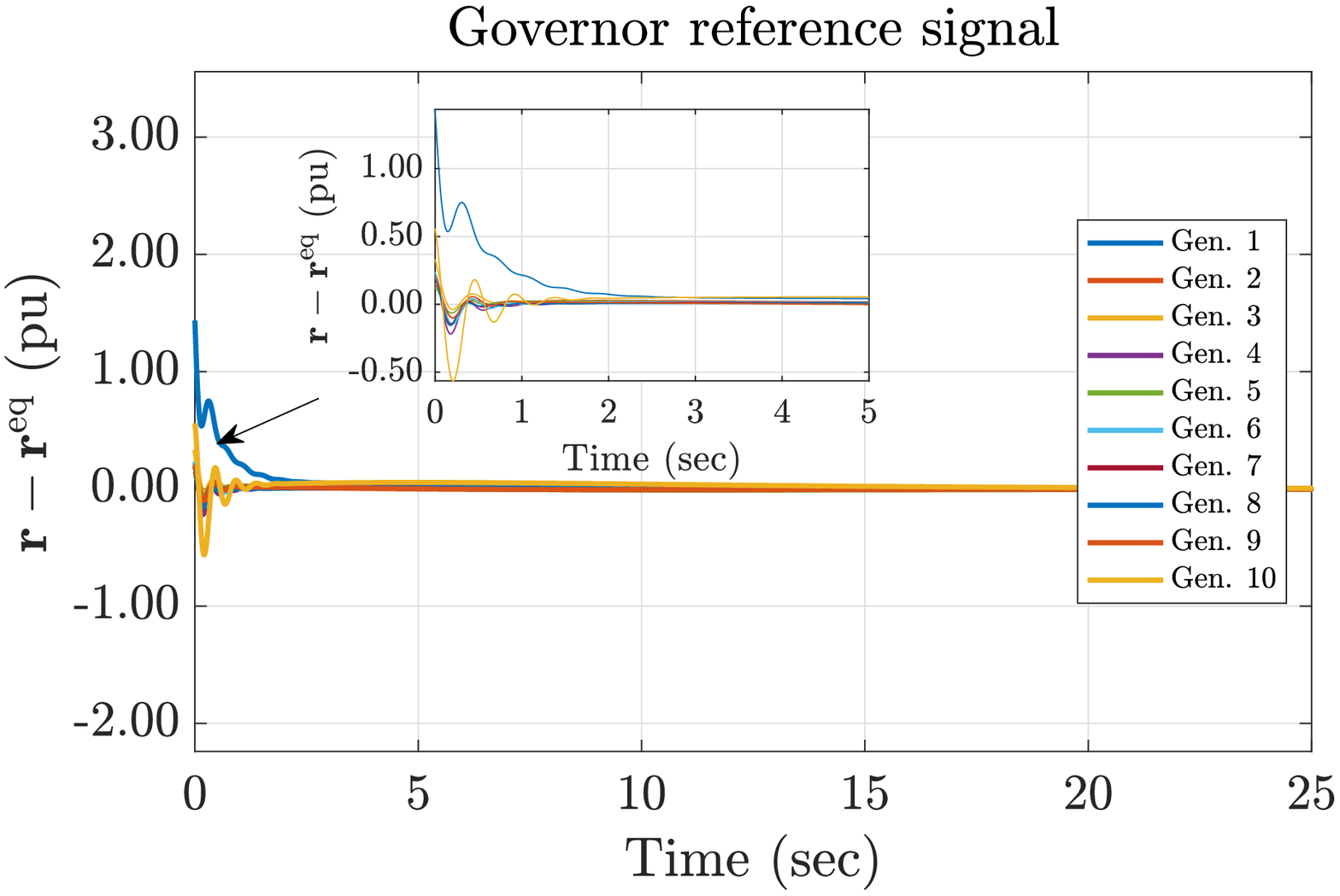}\label{fig:ALQR-OPFr}}{}  \quad
	\subfloat[]{\includegraphics[scale=0.3]{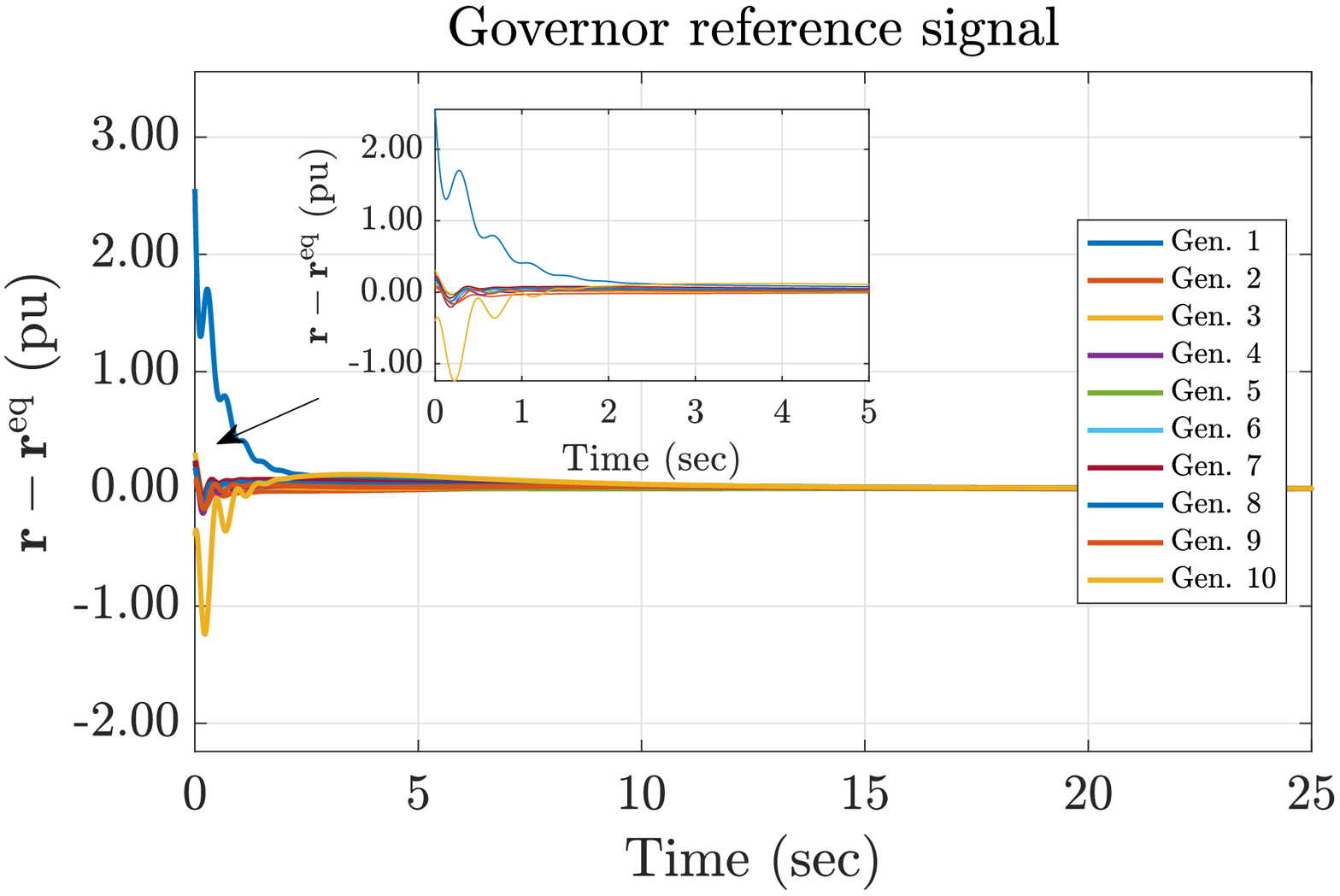}\label{fig:OPFr}}{}
	\caption{Generator frequencies and governor signals using LQR-OPF (left), ALQR-OPF (center), and OPF (right).} 	\label{fig:FreqMechGov}	
\end{figure*}

\begin{figure*}[!ht]
	\centering
		\subfloat[]{\includegraphics[scale=0.3]{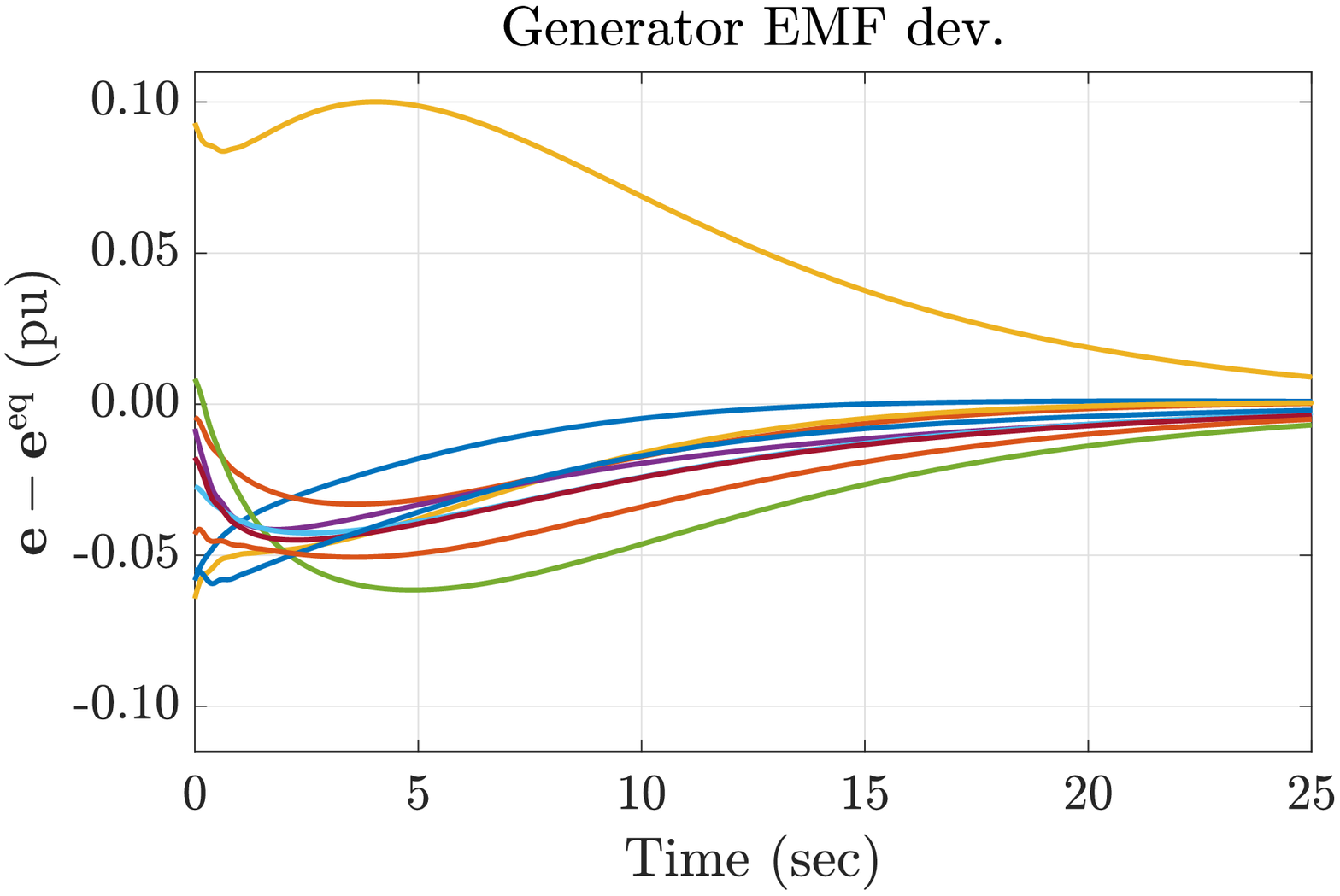}\label{fig:LQR-OPFe}}{} \quad
	\subfloat[]{\includegraphics[scale=0.3]{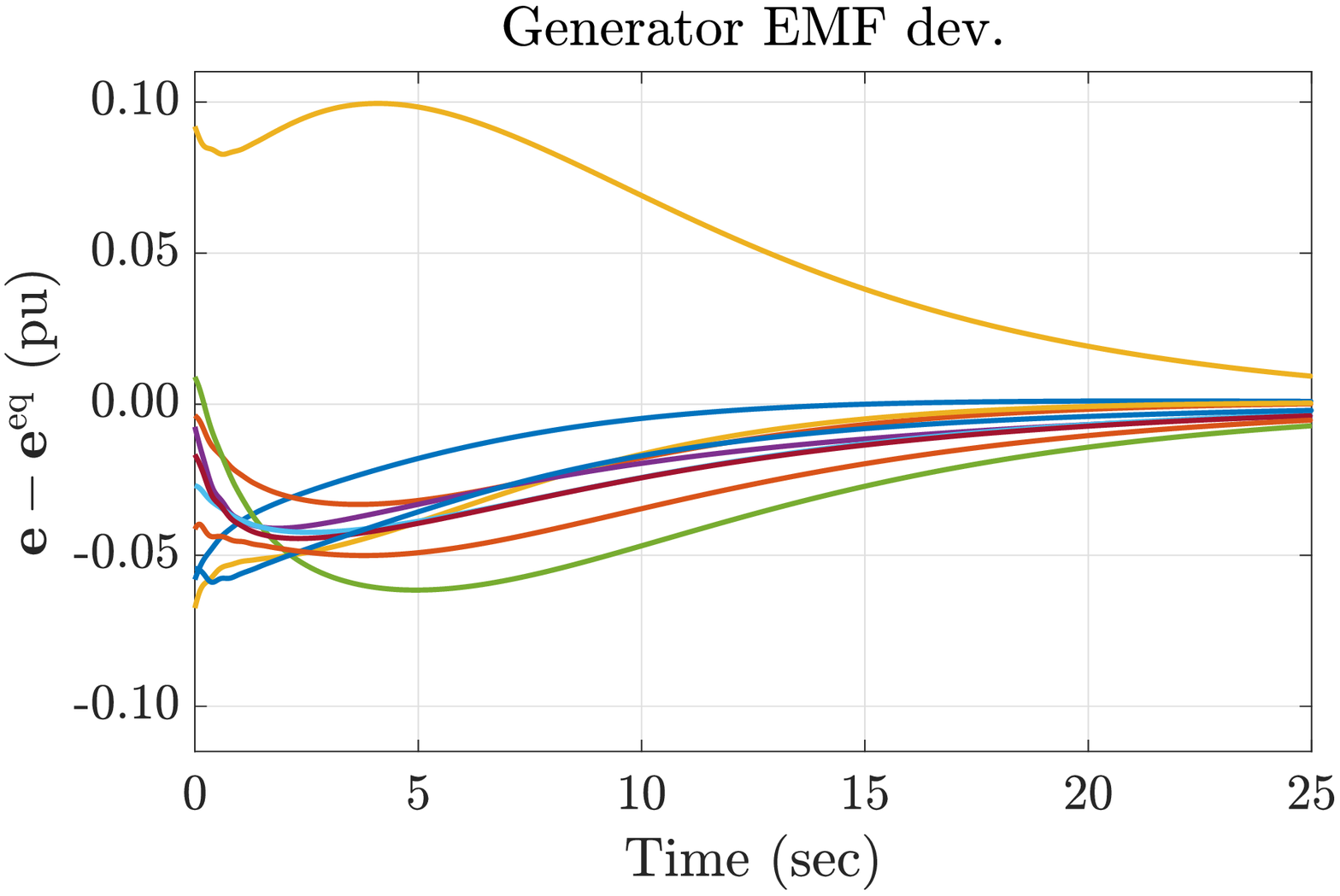}\label{fig:ALQR-OPFe}}{}  \quad
	\subfloat[]{\includegraphics[scale=0.3]{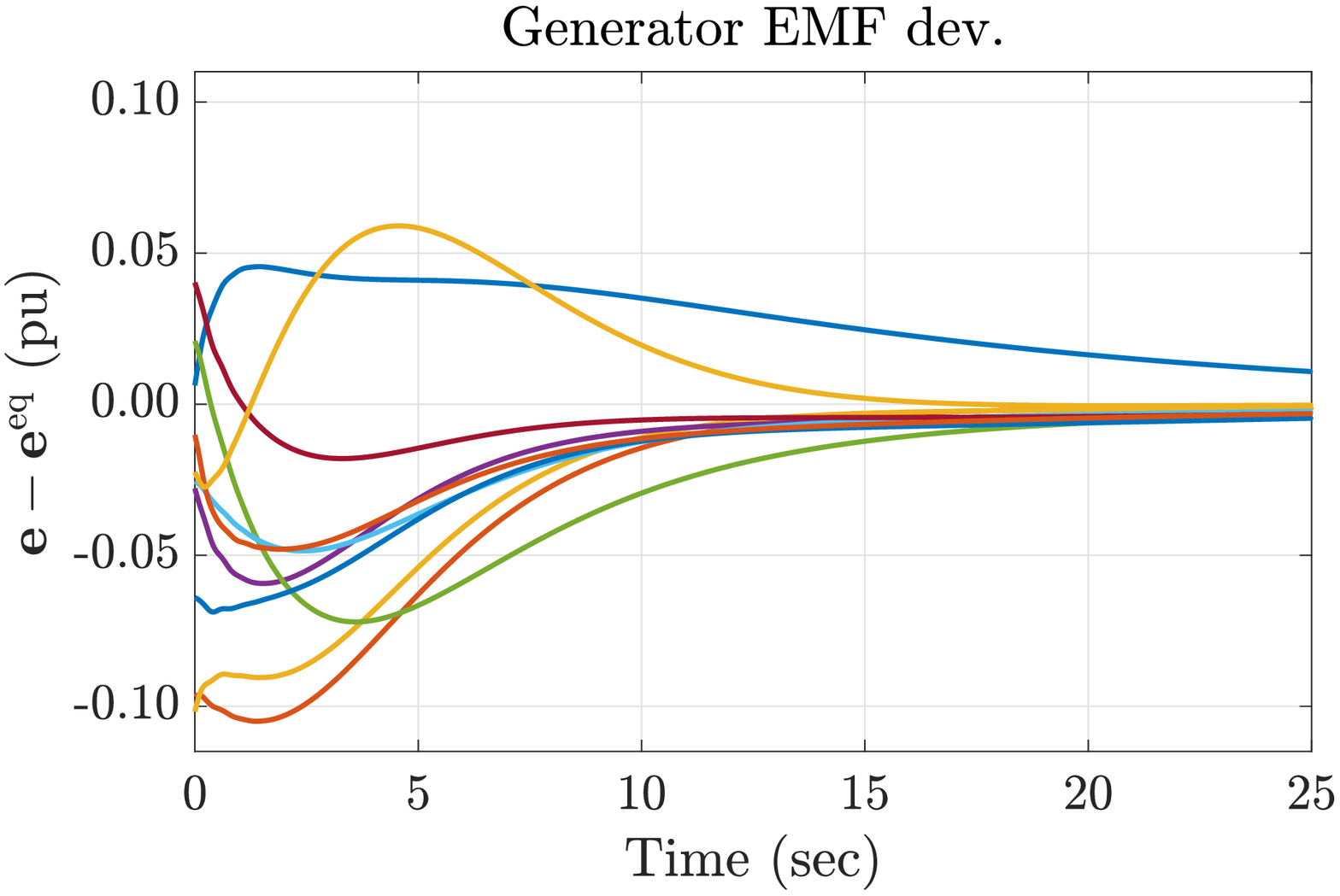}\label{fig:OPFe}}{} \\
		\subfloat[]{\includegraphics[scale=0.3]{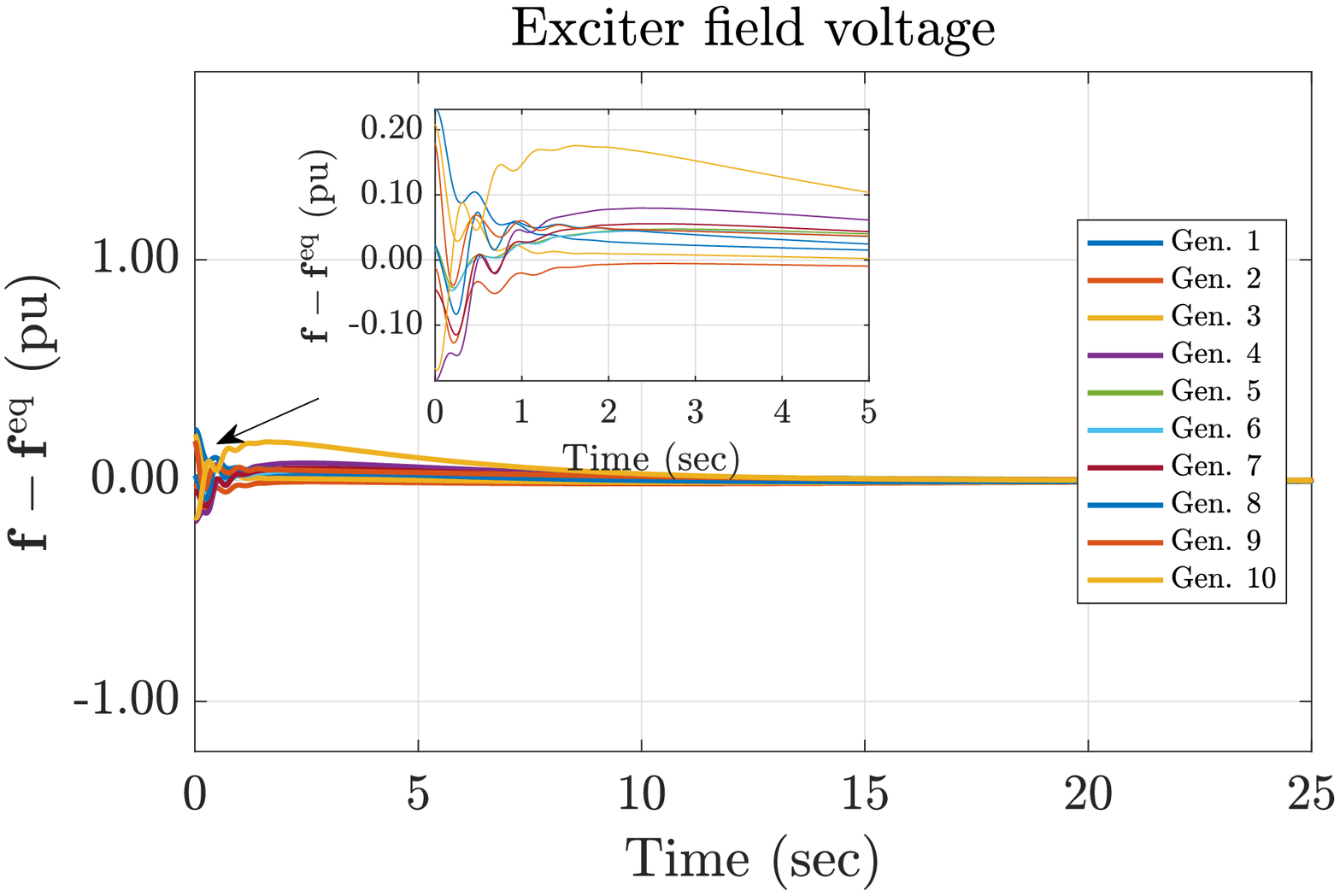}\label{fig:LQR-OPFf}}{}  \quad 
	\subfloat[]{\includegraphics[scale=0.3]{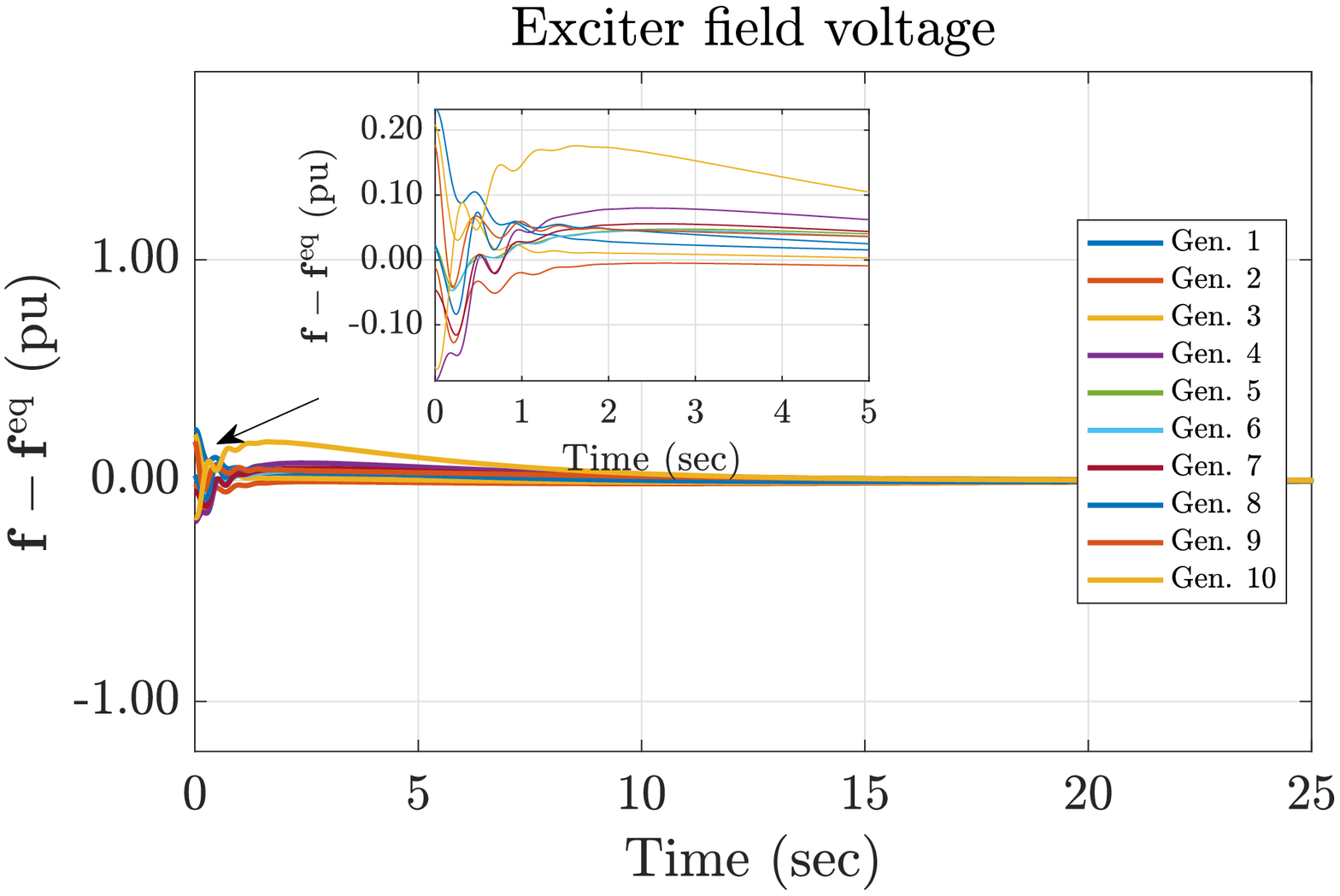}\label{fig:ALQR-OPFf}}{}  \quad
	\subfloat[]{\includegraphics[scale=0.3]{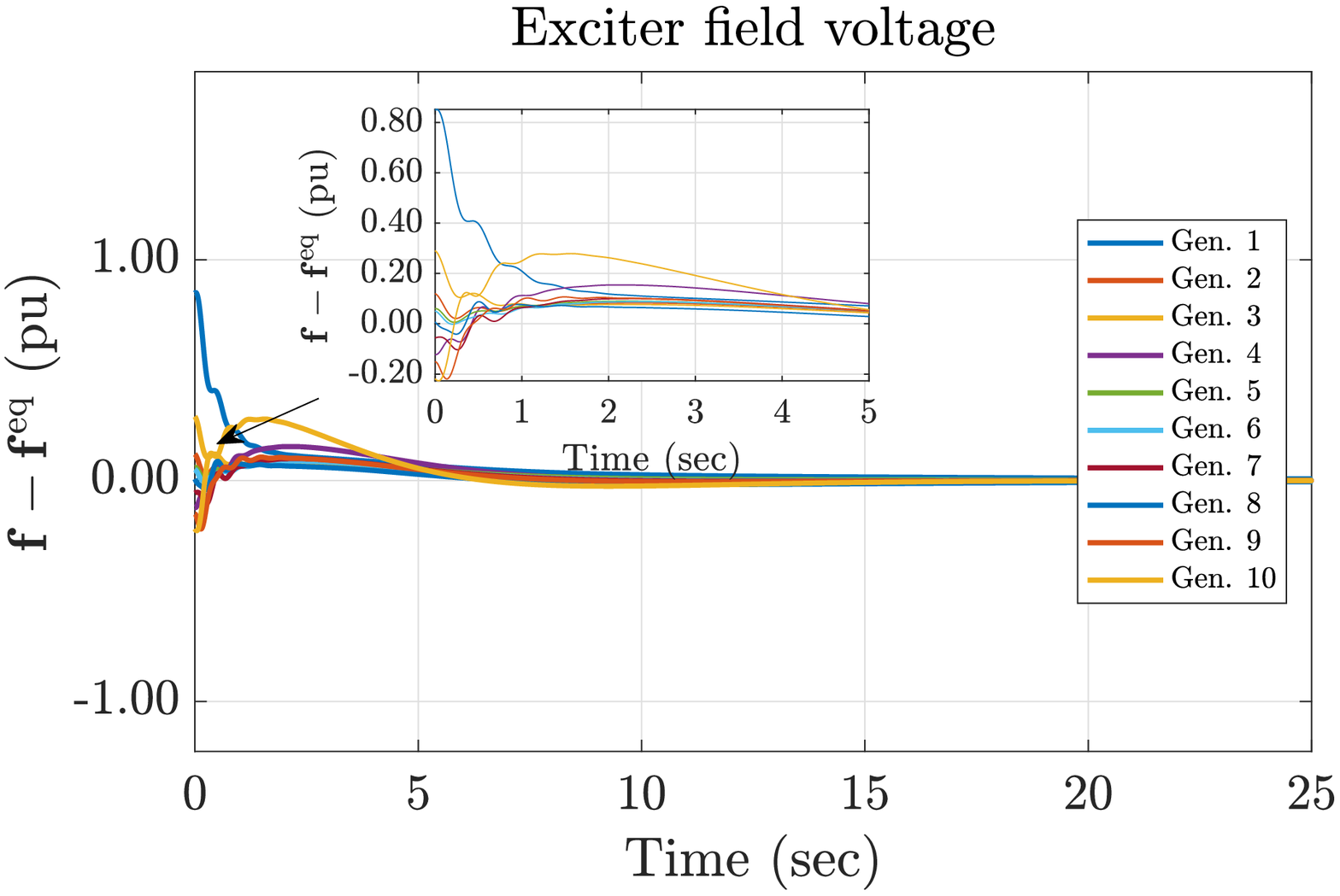}\label{fig:OPFf}}{}  \\
			\subfloat[]{\includegraphics[scale=0.3]{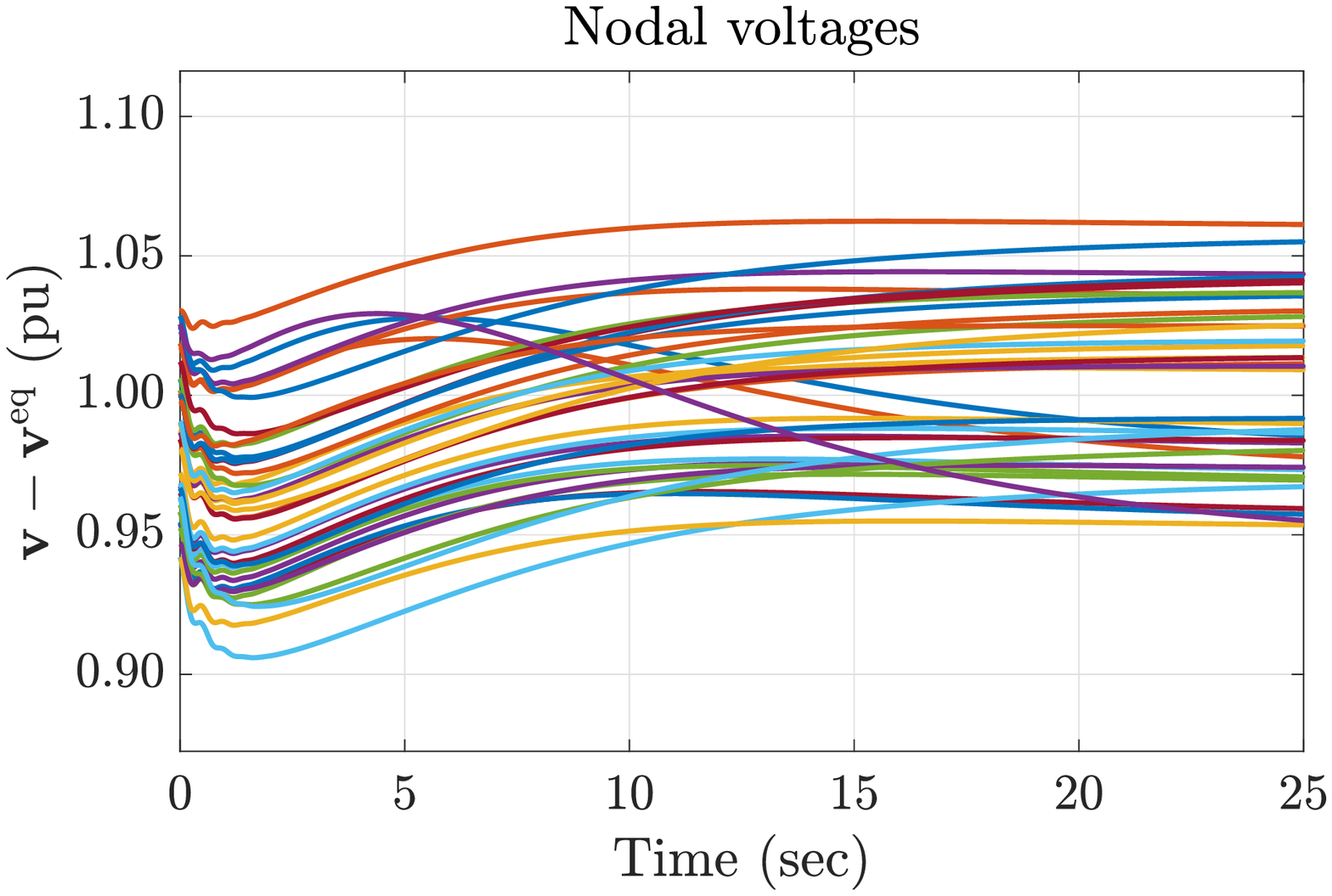}\label{fig:LQR-OPFv}}{}  \quad \:
	\subfloat[]{\includegraphics[scale=0.3]{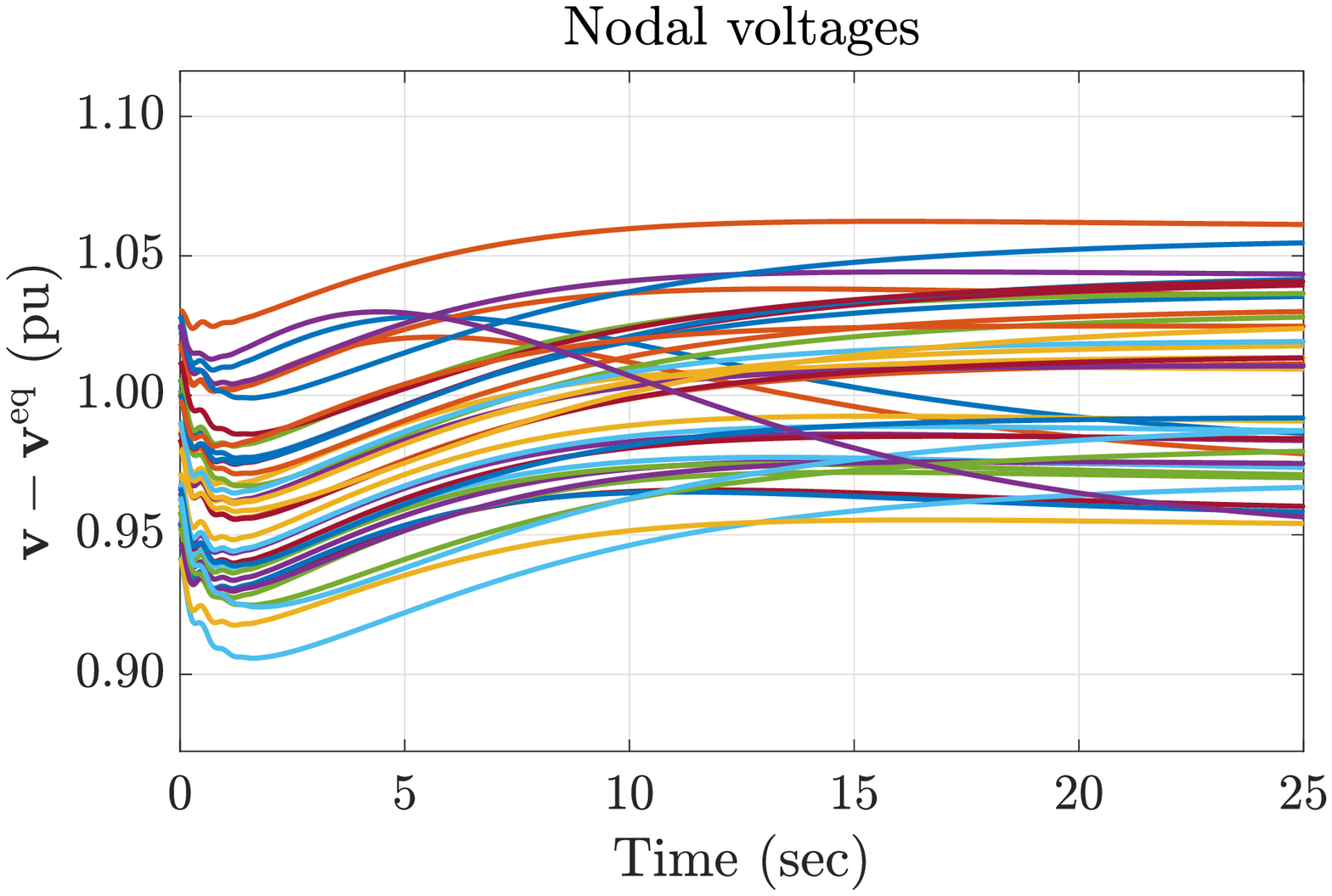}\label{fig:ALQR-OPFv}}{}  \quad 
	\subfloat[]{\includegraphics[scale=0.3]{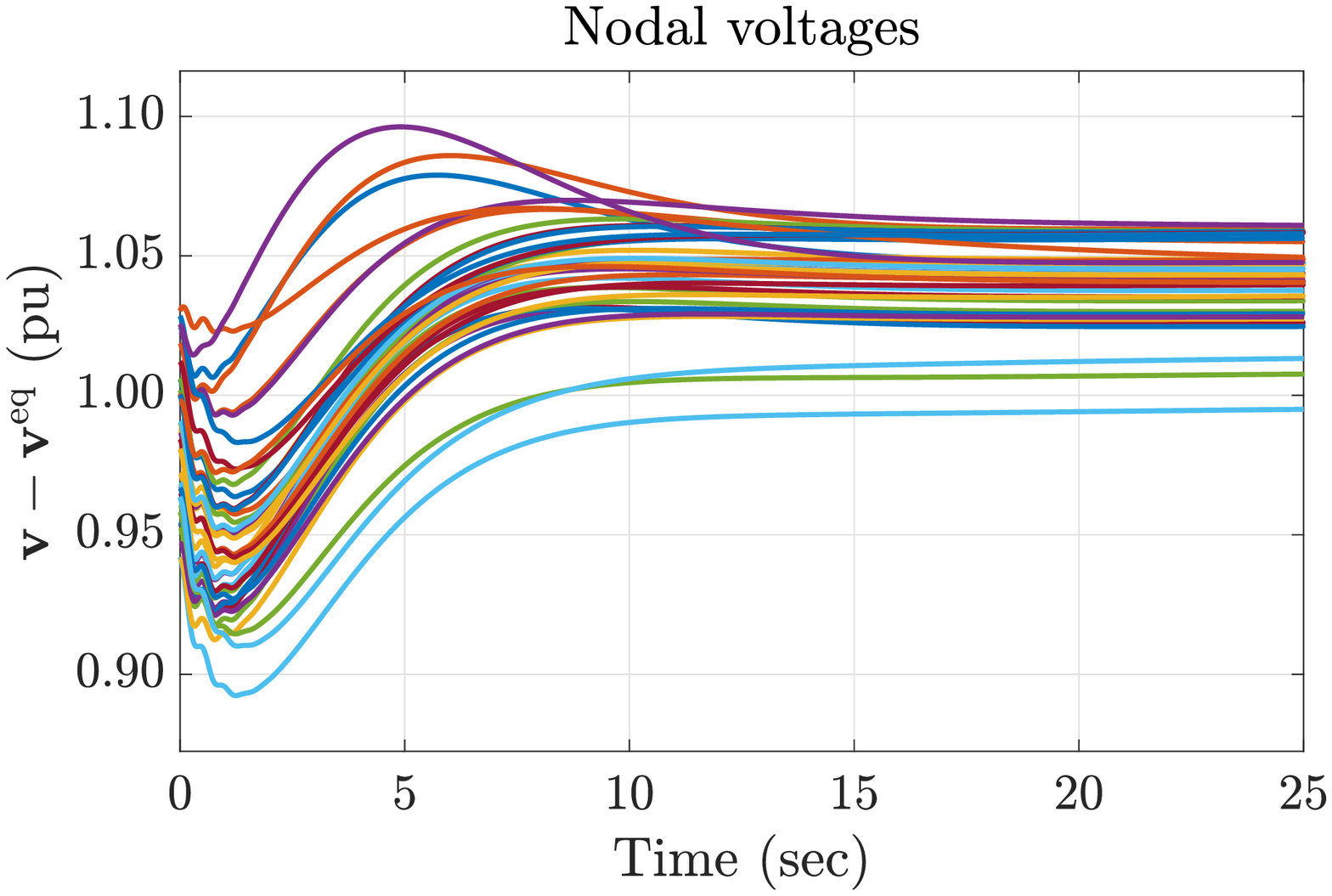}\label{fig:OPFv}}{} 
	\caption{ Generator EMF, exciter voltage,  and network nodal voltages using LQR-OPF (left), ALQR-OPF (center), and OPF (right).}
		\label{fig:EMFfVolt}	
\end{figure*}
\subsection{Dynamical simulation with LQR}
\label{sec:numtests:dynlqr}
A  step load increase of   $10\%$ in real power with power factor $0.9$ is applied at $t=0$. 
This implies that $\Delta p_{l_i}^s= 0.1 p_{l_i}^0$ and $\Delta q_{l_i}^s=0.0484 q_{l_i}^0$,  totaling  to a  load increase of $31.50 + j 5.56~\mr{MVA}$ for the 9-bus system, $25.90+j3.56~\mr{MVA}$ for the 14-bus system,  $625.42+j67.10~\mr{MVA}$ for the 39-bus system, $125.08+ j 16.28~\mr{MVA}$ for the 57-bus, $222.87 + j30.74~\mr{MVA}$ for the 200-bus Illinois system.
 This load increase drives the nonlinear dynamics \eqref{eqngroup:DAEs} out of the initial equilibrium. By applying LQR control according to the workflow in Fig.~\ref{fig:Diagram}, the dynamics in~\eqref{eqngroup:DAEs} are steered to arrive at the desired equilibria obtained from the OPF and LQR-OPF. 

With selections of $\alpha=0.6$ and $T_{\mr{lqr}}=1000$, Table \ref{table:costs} lists the  breakdown of steady-state, control, and total costs, as well as maximum frequency and voltage deviations from the optimal equilibrium of the LQR-OPF, ALQR-OPF with two iterations, and OPF. Column 3 gives the optimal objective of the LQR-OPF  in \eqref{eqngroup:lqropf} or the best objective found by the ALQR-OPF in Algorithm \ref{algorithm:ALQR-OPF}.  Observe that the optimal objectives of LQR-OPF and ALQR-OPF are almost identical, which implies that the approximation in Algorithm \ref{algorithm:ALQR-OPF} does not compromise optimality. 

Steady-state costs in column 4 correspond to $\mb{c}(\mb{a}^{\mr{eq}})$ found by each approach.  The LQR step in the diagram of Fig.~\ref{fig:Diagram} is solved using MATLAB's \texttt{care.m} by inputting $\mb{Q}(\mb{z}^{\mr{eq}})$ and $\mb{R}(\mb{z}^{\mr{eq}})$.  This process yields $\mb{P}$ and the corresponding feedback gain $\mb{K}=-\mb{R}^{-1}\mb{B}^{\top} \mb{P}$.  The estimates of the control costs are then calculated as $\frac{{T}_{\mr{lqr}}}{2}(\mb{x}^{\mr{eq}}- \mb{x}^0)^{\top} \mb{P}(\mb{x}^{\mr{eq}}- \mb{x}^0)^{\top}$  and are given in column 5 of Table \ref{table:costs}. The total estimated costs are then the summation of steady-state and estimated control costs and are given in column 6.  The computation times of LQR-OPF, ALQR-OPF, and OPF are listed in column 7. Notice that for the large 200-bus Illinois network, the LQR-OPF takes approximately 12 hours, while the ALQR-OPF solves the problem in less than three seconds and without significant loss in optimality. 

 Control costs reported in column 8 of Table~\ref{table:costs} are computed as $\frac{T_\mr{lqr}}{2}\int\nolimits_{0}^{t_f} (\Delta \mb {x'}^{\top} \mb{Q} \Delta \mb {x}' + \Delta \mb{u'}^{\top} \mb{R} \Delta \mb {u}')dt$, that is,  through numerical integration of the trajectories resulting from the simulation of the nonlinear DAEs.   The total cost, given in column 9, is simply the summation of control and steady-state cost. Between the coupled and decoupled approaches, OPF exhibits lower steady-state cost but higher control cost. In terms of total cost, the LQR-OPF and ALQR-OPF show improved performance.  The  maximum frequency deviation is also much lower for LQR-OPF and ALQR-OPF than the OPF.   

In Figs.~\ref{fig:FreqMechGov} and \ref{fig:EMFfVolt}	  the dynamical performance of the 39-bus system under LQR-OPF, ALQR-OPF, and OPF in conjunction with load-following LQR is depicted. Specifically, generator frequencies and governor reference signals are portrayed in Fig.~\ref{fig:FreqMechGov} where it is observed that quantities, especially frequencies, resulting from OPF undergo higher fluctuations than  those resulting   from LQR-OPF and ALQR-OPF.  In Fig.~\ref{fig:EMFfVolt} the generator internal EMF, exciter field voltage, and nodal voltages are depicted. Notice that overall, voltages obtained from LQR-OPF and ALQR-OPF exhibit smaller deviations from the designated equilibria in comparison with those obtained from OPF. 
Deviations of generator angles $\bm{\delta}$ also exhibit similar behavior. Plots for the remaining quantities ($\bm\delta$ and $\mb m$) and the corresponding plots for the remaining test networks (including the 200-bus Illinois network) are available online \cite{GithubCodes}.

\subsection{Dynamical simulation with AGC}
The LQR-OPF furnishes a steady-state operating point with desirable stability properties. After the steady-state operating point has been computed, one does not have to necessarily use LQR as a controller, but one could rather implement another dynamic control law such as AGC or PI-control~\cite{ZhangDominguezGarcia2016,gloverbook}. The purpose of this section is to examine the control cost to drive the system to the equilibrium computed by LQR-OPF or OPF when the control law is AGC. 

To this end, the setup of Section \ref{sec:numtests:dynlqr} is followed, but AGC is used
 to adjust the governor reference signal during the dynamical simulations instead of LQR.  Dynamical equations describing the AGC for a multi-area power network are adopted from~\cite{ZhangDominguezGarcia2016} and~\cite{gloverbook}.   The selection of participation factors to steer the DAEs to the desired equilibrium follows the suggestions in \cite[p.~87]{wollenbergbook2012} for computer implementations. The specifics are detailed in the Appendix.

\begin{table}[!ht]
	\scriptsize
	\centering
		\renewcommand{\arraystretch}{1.2}
		\caption{Costs Comparison between LQR-OPF, ALQR-OPF, and OPF methods at $\alpha=0.6$ under AGC }
		\label{table:AGCcosts}
		\begin{tabular}{|c|c|c|c|c|c|}
			\hline
			\multirow{3}{*}{Network}&	\multirow{3}{*}{Method}      &     Control & Total      & Max.&  Max. \\
			&  & cost &  cost & freq. dev.  &  volt. dev.  \\
			&   & ($\$$)   & ($\$$)  & ($\mr{Hz}$) &   ($\mr{pu}$)  \\
			\hline
			\hline
				\multirow{3}{*}{9-bus}    & LQR-OPF                   & 35.92           & 6180.04         & 0.0144          & 0.0262    \\       
  & ALQR-OPF               & 35.77           & 6180.02         & 0.0144          & 0.0261          \\ 
		   & OPF                    & 583.90          & 6697.50         & 0.0161          & 0.1240     \\
			\hline 
			\multirow{3}{*}{14-bus}  & LQR-OPF                & 94.19           & 9273.01         & 0.0040          & 0.0297          \\ 
		 & ALQR-OPF                & 94.30           & 9272.86         & 0.0041          & 0.0297       \\    
			  & OPF                       & 590.84          & 9718.19         & 0.0096          & 0.0297      \\
			\hline
			\multirow{3}{*}{39-bus}   & LQR-OPF                 & 4819.03         & 57691.30        & 0.0765          & 0.0782    \\       
  & ALQR-OPF                & 4843.82         & 57715.62        & 0.0765          & 0.0772       \\    
		 & OPF                     & 17781.72        & 69167.74        & 0.1321          & 0.1062    \\
			\hline
			\multirow{3}{*}{57-bus}         & LQR-OPF                & 4507.44         & 52829.71        & 0.0489          & 0.0599    \\       
	       & ALQR-OPF                & 4401.21         & 52769.55        & 0.0481          & 0.0600   \\        
	         & OPF                    & 12652.76        & 59852.51        & 0.0760          & 0.0637 \\
			\hline
			\multirow{3}{*}{200-bus}  & LQR-OPF        & 9214.91         & 57561.95        & 0.0249          & 0.0608           
			\\
			& ALQR-OPF               & 5667.92         & 56017.57        & 0.0168          & 0.0589     \\      
			 & OPF                & 14959.59        & 63231.30        & 0.0358          & 0.0741    \\ 
			\hline
	\end{tabular}
\end{table}

\begin{figure}[t]
	\centering
\includegraphics[scale=0.43]{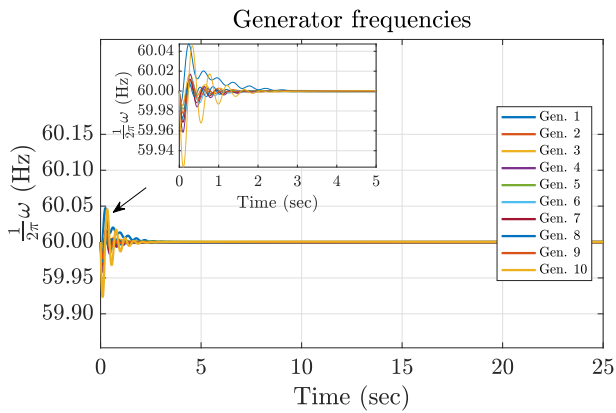}
	\caption{Generator frequencies using  ALQR-OPF in conjunction with AGC.}
	\label{fig:ALQR-OPF-AGC-freq}
\end{figure}

\begin{figure}[t]
\centering
	\includegraphics[scale=0.43]{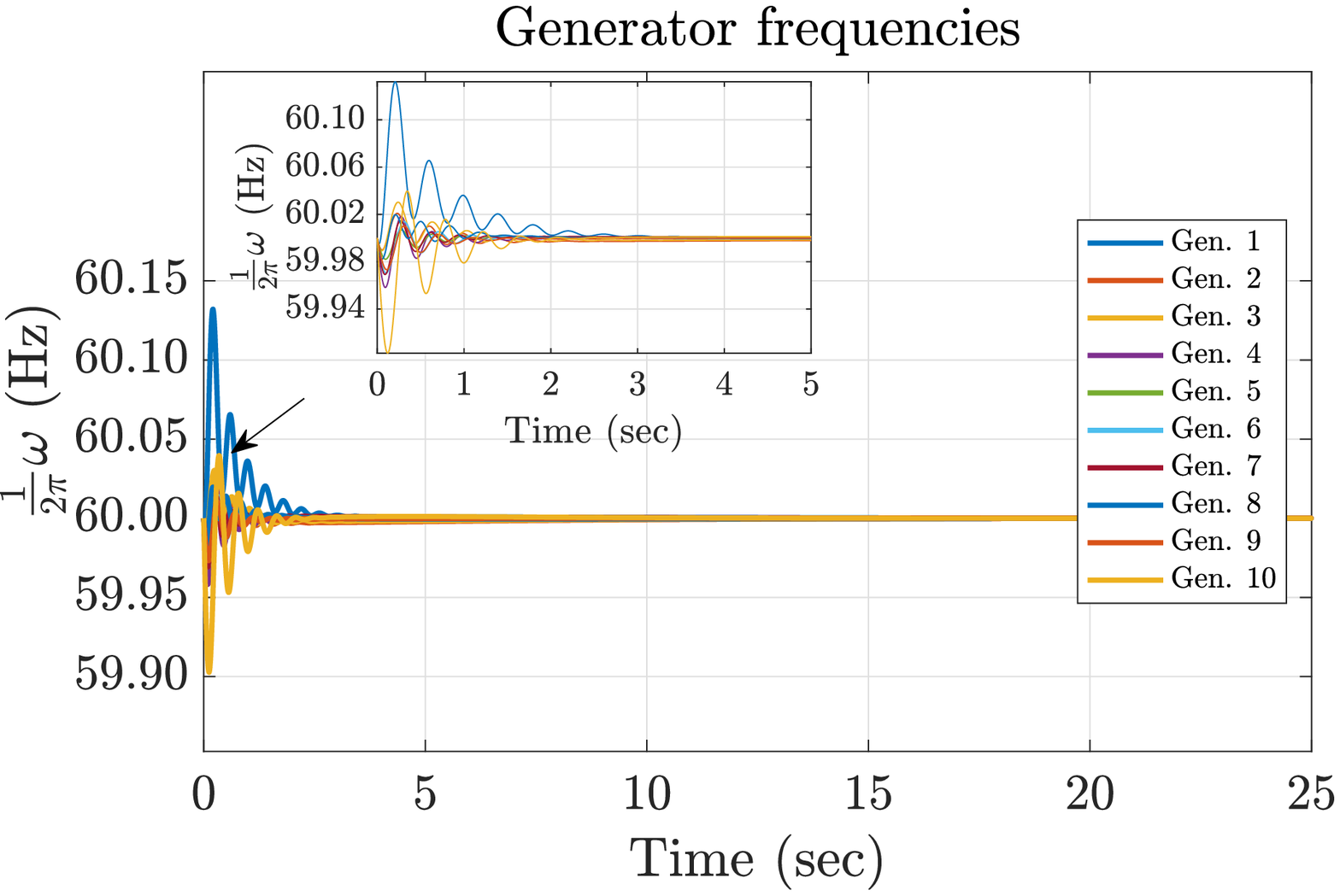}
	\caption{Generator frequencies using  OPF in conjunction with AGC.}
	\label{fig:OPF-AGC-freq}
\end{figure}

\begin{figure}[t]
\centering
\includegraphics[scale=0.43]{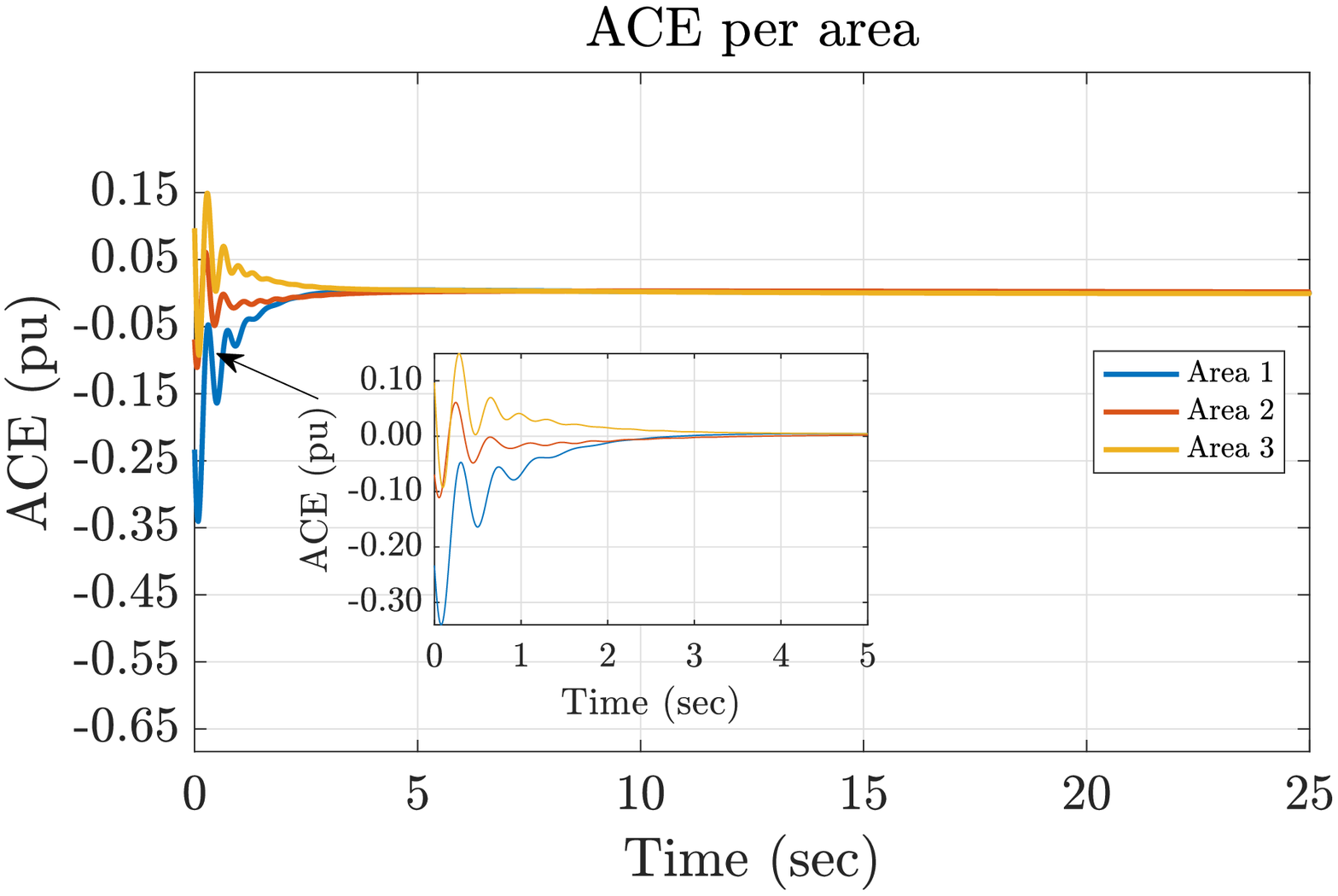}
\caption{Area control signal using  ALQR-OPF in conjunction with AGC. }
\label{fig:ALQR-OPF-AGC-ACE}
\end{figure}

\begin{figure}[t]
\centering
	\includegraphics[scale=0.43]{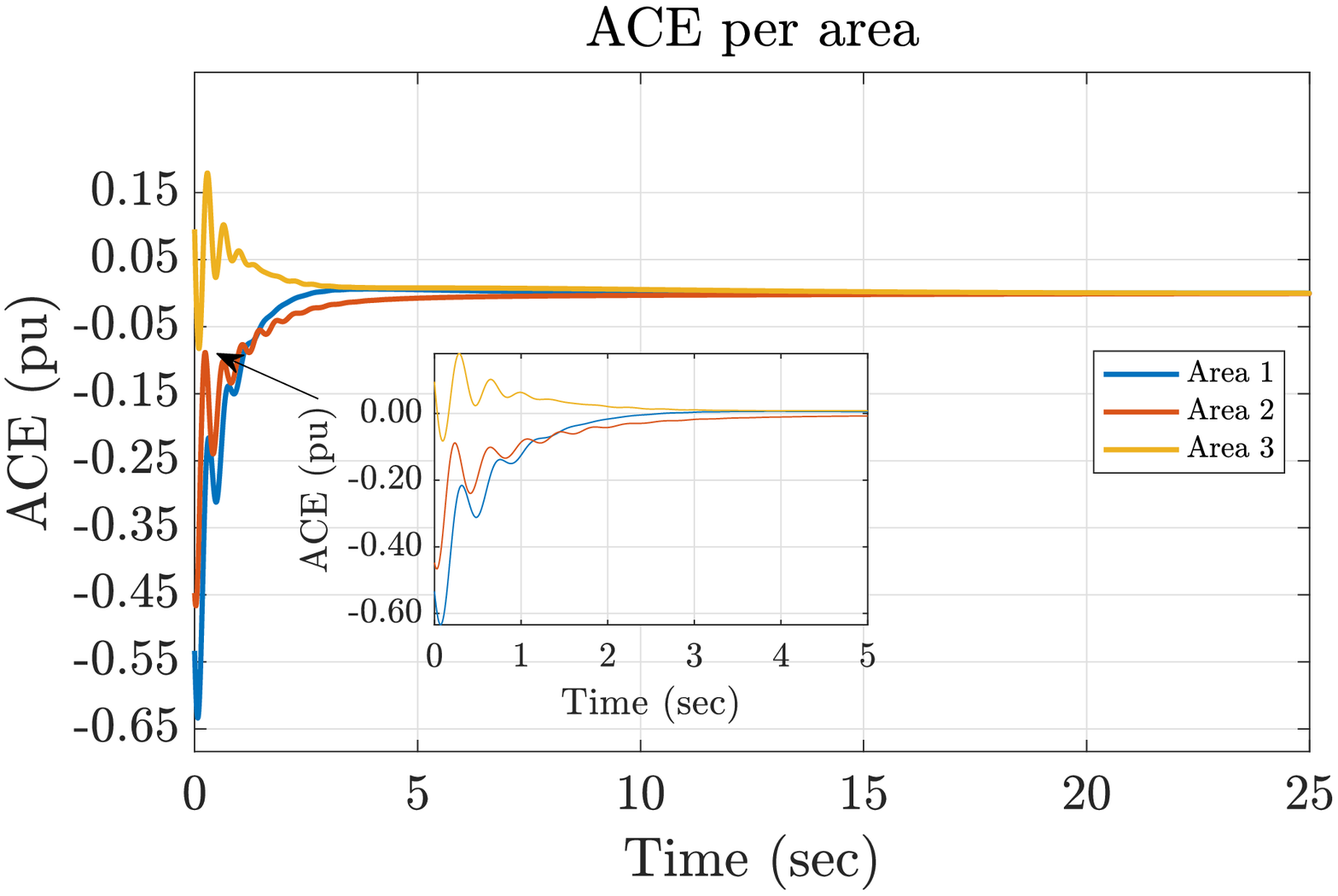}
	\caption{Area control signal using OPF in conjunction with AGC. }
	\label{fig:OPF-AGC-ACE}
	\end{figure}

Table \ref{table:AGCcosts} reports the control and total costs of LQR-OPF and ALQR-OPF compared with OPF, when AGC is used to adjust the governor signals $r_i$ for $i \in \mc{G}$. Notice that the steady-state results are those previously given in Table \ref{table:costs} and only quantities pertaining to dynamical simulations have changed. It is observed that using AGC, setpoints provided by LQR-OPF and ALQR-OPF  result in smaller control costs compared to OPF.\footnote{Table \ref{table:AGCcosts} indicates that total costs of ALQR-OPF are sometimes smaller than those of LQR-OPF. This happens because the optimality of LQR-OPF over ALQR-OPF is only guaranteed when the associated controller is LQR. 
}   The corresponding dynamical performance for ALQR-OPF and OPF in conjunction with AGC is provided in Fig.~\ref{fig:ALQR-OPF-AGC-freq}--\ref{fig:OPF-AGC-ACE}. The setpoints provided by ALQR-OPF result in smaller frequency deviation and smaller ACE signal compared to the setpoints provided by OPF. The performance using the setpoints of LQR-OPF is similar to that of ALQR-OPF and has been omitted for brevity. Corresponding plots of other system quantities for the 39-bus network and the remaining test networks are also made available online~\cite{GithubCodes}.

\color{black}
\subsection{Effect of coupling}
\begin{figure}[t]
	\centering
	\includegraphics[scale=0.45]{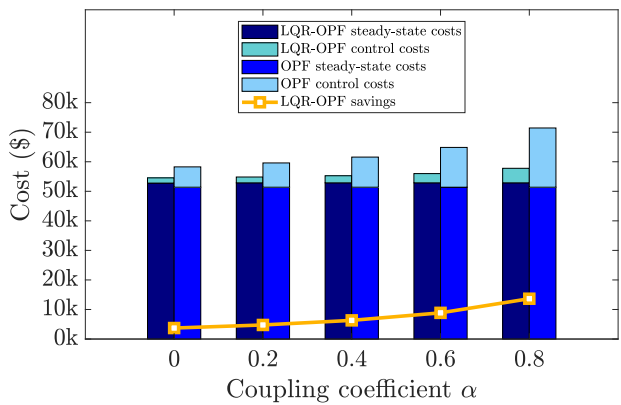}
	\caption{Effect of coupling between steady-state variables and control cost matrices $\mb{Q}$ and $\mb{R}$.   As coupling coefficient $\alpha$ increases, LQR-OPF results in  more savings than the decoupled OPF and LQR approach. }
	\label{fig:coupling}	
\end{figure}

Here, the effect of parameter $\alpha$ which couples  steady-state variables to load-following control costs through \eqref{eqngroup:QRinv} is studied. When the value of $\alpha$ increases to approach $1$, entries of matrices $\mb{Q}$ and $\mb{R}$ increase as the values of $p_{g_i}^s$ and $q_{g_i}^s$ approach their respective maximum. It is depicted in Fig.~\ref{fig:coupling} that as the coupling coefficient $\alpha$ increases, control costs increase in both OPF and LQR-OPF. However, control costs of LQR-OPF are significantly lower than the costs incurred by the scheme where OPF and LQR are
solved independently.

\subsection{ALQR-OPF on larger networks}
This section examines the performance of ALQR-OPF on larger networks. Table~\ref{table:LargerNetworks} lists the steady-state costs, estimated control costs (as computed by $\frac{T_{\mr{lqr}}}{2}(\mb{x}^{\mr{eq}}- \mb{x}^{0})^{\mr{\top}} \mb{P} (\mb{x}^{\mr{eq}}- \mb{x}^0)$), total estimated costs, and computation times. The main observation is that ALQR-OPF yields significantly smaller total costs than OPF. The computation time of ALQR-OPF is larger than that of OPF, but it is worth noting that ALQR-OPF has been solved by a general-purpose solver, while MATPOWER's solver is specifically tailored to the OPF problem.

\begin{table}[t]
	\centering
	\scriptsize
	\caption{ALQR-OPF on larger networks}
	\label{table:LargerNetworks}
	\begin{tabular}{|c|c|c|c|c|c|}
		\hline
	\multirow{3}{*}{Network}&	\multirow{3}{*}{Method}      & Steady-state &    Control & Total   & Comp.     \\
& & cost & est. cost &  est. cost & time \\
& &  ($\$$) & ($\$$)   & ($\$$) & (seconds)   \\
\hline
\hline
\multirow{2}{*}{1354-bus}   & ALQR-OPF           & 82612        & 2311022      & 2393635      & 72.65   \\
 & OPF             & 81688        & 8279013      & 8360702      & 6.06       \\      
\hline
\multirow{2}{*}{2383-bus}  & ALQR-OPF              & 2235711     & 21174       & 2256885      & 325.72    \\
   & OPF               & 2217287     & 123999       & 2341287      & 9.86     \\         
\hline
\multirow{2}{*}{2869-bus} &  ALQR-OPF               & 149453       & 3255044      & 3404498     & 683.74           \\
 & OPF              & 147865       & 14830044     & 14977910     & 52.70      \\      
\hline
	\end{tabular}
\end{table}

\section{Summary and Future Work}
\label{sec:conclusion}
An OPF framework is presented that in addition to solving for optimal steady-state setpoints   provides an optimal feedback law to perform load-following control.  The costs of load-following control  is captured by a classical LQR control that accounts for  deviations of system states and controls from their optimal steady-state setpoints.  A joint formulation  of OPF and load-following control, termed LQR-OPF, is obtained by combining a linearized OPF with an equivalent  SDP  formulation of the LQR.   Numerical tests verify that compared to a scheme where OPF and load-following control problems are solved separately,  LQR-OPF  features significantly improved dynamic performance and reduced overall system costs.  

The proposed framework  is general and allows for seamless incorporation of different power system applications, such as wind turbine and storage device dynamics---both operating at different time-scales.  Recent work,  for instance, demonstrates the impact of wind power injection on power system oscillations~\cite{ChandraGaymeChakrabortty2016}. Future work includes integrating more modern  applications into the proposed framework while investigating the regulation and cost benefits of this integration.

\appendix[AGC implementation]
Denote  by $\mc{A}$  the set of areas of a power  network and by  $\mc{A}_a$  the set of neighboring areas to area $a \in \mc{A}$. Further, denote respectively by $p_{aa'}$ and $p_{aa'}^{\mr{eq}}$, the time-varying and equilibrium  aggregate real power flows  on the tie-lines from  area $a \in \mc{A}$ into area $a' \in \mc{A}$ . Notice that $p_{aa'}$ is a function of the voltage magnitudes and angles at the terminals of all the tie-lines connecting $a$ to $a'$. 
The area control error (ACE) for area $a \in \mc{A}$ is then computed as
\begin{IEEEeqnarray}{rCl}
	\mr{ACE}_a=\sum\limits_{a' \in \mc{A}_a} (p_{aa'}-p_{aa'}^{\mr{eq}})+b_a \left(\frac{1}{|\mc{G}_a|} \sum\limits_{i \in \mc{G}_a} \omega_i - \omega^{\mr{s}}\right)  \IEEEeqnarraynumspace \label{eqn:ACE}
\end{IEEEeqnarray}
where $\mc{G}_a$ denotes the set of generators in area $a$, and  $b_a=\sum_{i \in \mc{G}_a} (\frac{1}{R_i}+D_i)$ is the area bias factor.

AGC uses the ACE signal to provide a command to the governor reference signal. The corresponding dynamical equations for area $a \in \mc{A}$ are 
\begin{subequations} 
	\label{eqngroup:AGC}
	\begin{IEEEeqnarray}{rCl}
		\dot{y}_a &=& K_{a} \left( -y_{a} -  \mr{ACE}_{a}   + \sum\limits_{i \in \mc{G}_a} p_{g_i}^{\mr{eq}}\right) \label{eqn:ydotAGC} \\
		r_i&=&K_i y_a,  \quad  i \in \mc{G}_a, \label{eqn:rAGC}
	\end{IEEEeqnarray}
\end{subequations}
where $K_a$ is an integrator gain, $K_i$ is the participation factors of generator $i$, and $r_i$ is fed back into~\eqref{eqn:mdynamics}. The participation factors are set to $K_i = p_{g_i}^{\mr{eq}}/ \sum_{i \in \mc{G}_a} p_{g_i}^{\mr{eq}}$ so that the power network will be steered to the desired equilibrium.  Notice that the sum of participation factors per area equals unity.
\bibliographystyle{IEEEtran}
\bibliography{references}

\end{document}